\documentclass[12pt]{article}
\usepackage{amsmath,amssymb,amscd}
\usepackage{bm}

\usepackage{amsmath}
\usepackage{graphicx}

\newtheorem{tm}{Theorem}[section]
\newtheorem{lm}[tm]{Lemma}
\newtheorem{co}[tm]{Corollary}
\newtheorem{re}[tm]{Remark}

\newtheorem{pr}[tm]{Proposition}
 
 \newenvironment{demo}[1]{\par\smallskip\par\begin{trivlist}
\item[]{\bf #1}\ }{\end{trivlist}\par\smallskip\par}
\newcommand{\Proof}{\begin{demo}{{\it Proof.\ }}}
\newcommand{\qed}{\end{demo}}
\newcommand{\toy}{\ \rule[0em]{0.5ex}{1.8ex}}
\newcommand{\QED}{\toy\end{demo}}
\newcommand{\la}{\langle}
\newcommand{\ra}{\rangle}
\newcommand{\nn}{\nonumber}
\newcommand{\III}{{\vert \kern-.10em \vert \kern-.10em \vert}}
\newcommand{\ve}{\varepsilon}
\newcommand{\al}{\alpha}

\makeatletter
 
 \@addtoreset{equation}{section}
\makeatother
 
\begin{document}
\setlength{\baselineskip}{15pt} 
%%%%%%%%%%%%%%%%%%%%%%%%%%%%%%%%%%%%%%%%%%%%%%%%%%%%%%%%%%%%%%%%%%%%%%
%
\bibliographystyle{plain}
\title{
Large deviations for rough path lifts of
 Watanabe's pullbacks of delta functions 
\footnote{
%Preliminary version: ~ 26 Dec 2014.
%
%
{\bf Mathematics Subject Classification:}~ 60F10, 60H07, 60H99, 60J60.
{\bf Keywords:} large deviation principle, rough path theory, Malliavin calculus,
quasi-sure analysis,
pinned diffusion process.
}
}
%%%%%%%%%%%%%%%%%%%%%%%%%%%%%%%%%%%%%%%%%%%%%%%%%%%%%%%%%%%%%%%%%%
%%%%%%%%%%%%%%%%%%%%%%%%%%%%%%%%%%%%%%%%%%%%%%%%%%%%%%%%%%%%%%%%%%%
\author{ 
Yuzuru Inahama
%\footnote{
%Preliminary version: ~ 10 Dec 2014.
%
%e-mail:~\tt{inahama@math.nagoya-u.ac.jp} 
%\\
%Graduate School of Mathematics,   Nagoya University
%\\
%Furocho, Chikusa-ku, Nagoya 464-8602, JAPAN.
%\\
%E-mail:~\tt{inahama@math.nagoya-u.ac.jp}
%}
%}
}
%\date{ \today }
\date{   }
% \pagestyle{empty}
%
% Start !!!
%
\maketitle
% \thispagestyle{empty}
%%%%\hspace{-5.5mm}

%%%%%%%%%%%%%%%%%%%%%%%%%%%%%%%%%%%%%%%%%%%%%%%%%%%%%%%%
%%%%%%%%%%%%%%%%%%%%%%%%%%%%%%%

\begin{abstract}
We study Donsker-Watanabe's delta functions 
associated with strongly hypoelliptic diffusion processes indexed by a small parameter. 
They are  finite Borel measures on the Wiener space 
and admit a rough path lift.  
Our main result  is  
a large deviation principle 
of Schilder type for the lifted measures
on the geometric rough path space
as the scale parameter tends to zero.
As a corollary, we obtain a large deviation principle 
conjectured by Takanobu and Watanabe, which is a generalization 
of a large deviation principle 
of Freidlin-Wentzell type for pinned diffusion processes.
\end{abstract}

\section{Introduction}

In 1993 Takanobu and Watanabe \cite{tw} presented a large deviation principle (LDP)
of Freidlin-Wentzell type 
for solutions of stochastic differential equations
(SDEs) under the strong  H\"ormander condition anywhere.
Unlike in the usual LDP of this type, 
the probability measures in \cite{tw} are not the push-forwards of the (scaled) Wiener measure, 
but the push-forwards of the measures of finite energy which is defined by 
the composition of the solutions of SDEs and the delta functions 
(i.e., Watanabe's pullbacks of the delta functions, also known as Donsker's delta function).
One interpretation of this LDP is a generalization of 
the LDP  of Freidlin-Wentzell type for pinned diffusion measures.
This LDP (Theorem 2.1, \cite{tw}) looks very nice. 
To the author's knowledge, however, no proof has been given yet.

% \vspace{5mm}

In this paper we reformulate this LDP on the geometric rough path space 
by lifting these measures in the rough path sense
and prove it rigorously by using quasi-sure analysis 
(which is a kind of potential theory in Malliavin calculus).  
Then, Theorem 2.1, \cite{tw} is a simple corollary  of  our main result.
After suitably specializing it, we also obtain 
the LDP for pinned diffusion measures under the strong H\"ormander condition anywhere.
Our main tools are rough path theory, Watanabe's  distributional Malliavin calculus,
and quasi-sure analysis.
%
%(Even this one seems new.)

%\vspace{5mm}

The elliptic case was already done in the author's previous work \cite{in2}.
This work is a generalization of it
to the strongly hypoelliptic case.
Note that 
many basic results on quasi-sure analysis for Brownian rough path were already obtained in \cite{in2}. 
%
%Our main tools are rough path theory, Watanabe's  distributional Malliavin calculus,
%and quasi-sure analysis.
%
Compared to \cite{in2},
the lower estimate becomes more difficult,  while the upper estimate remains somewhat similar.

Let us briefly recall the history of LDP of Schilder type 
on rough path space.
The first result was for the law of (scaled) Brownian rough path by
Ledoux, Qian, and Zhang \cite{lqz}.
By the continuity of Lyons-It\^o map in the rough path setting,
the usual LDP of Freidlin-Wentzell type is immediate.
Although a few nice proofs of this LDP 
are known, 
this new proof is quite straight forward and looks powerful.
Since then, LDPs of Schilder type became one of the central topics 
in the probabilistic aspects of rough path theory
and many papers have been written on it.
(For example,  an LDP for a wide class of Gaussian rough paths 
is proved in Friz and Victoir \cite{fv2}.
This class includes fractional Brownian rough path with Hurst parameter $H \in (1/4, 1/2]$.
The original motivation of \cite{in2} was to extend the idea in \cite{lqz} 
to the case of pinned diffusion measures.)
Another advantage of this method is that
one can also prove Laplace approximation (i.e., the precise asymptotics of LDP of Freidlin-Wentzell type)
along the same streamline with or without Malliavin calculus. 
(For example, see \cite{in3, in4} for results for fractional Brownian rough path).
In short, LDP theory on rough path space turned out to be quite successful.
Therefore, we believe that the geometric rough path space is the right place 
for the LDP conjectured by Takanobu and Watanabe in \cite{tw}.

%%%%%%%%%%%%%%%%%%%%%%%%%%%%%%
%%%%%%%%%%%%%%%%%%%%%%%%%%%%%%
%%
%\vspace{5mm}
%%%%%%%%%%%%%%%%%%%%%%%%%%%%%%
%%%%%%%%%%%%%%%%%%%%%%%%%%%%%%

The organization of this paper is as follows.
In Section 2 we introduce the setting, make basic assumptions, 
and state our main result (Theorem \ref{tm.main}) and its corollaries (Corollaries \ref{co.TW} and \ref{co.pFW}).
Section 3 is devoted to calculations of the skeleton ODE. 
We prove that
the deterministic Malliavin covariance matrix is non-degenerate 
at sufficiently many Cameron-Martin paths. 
This is a key lemma in the proof of the lower estimate
of our main theorem.
In Section 4 
we present some preliminaries on quasi-sure analysis on rough path space,
all of which were already shown or used in \cite{in2}.

We prove the large deviation lower estimate in Section 5.
Compared to the elliptic case in \cite{in2},
this part becomes more difficult for two reasons.
(These are closely related, however.)
One is non-degeneracy of 
the deterministic Malliavin covariance matrix. 
It fails at some Cameron-Martin paths in the hypoelliptic case.
(The aim of Section 3 is to deal with this difficulty).
The other is that 
uniform non-degeneracy of Malliavin covariance matrix 
of the diffusion processes does not hold in general.
We will use a modified version of the 
asymptotic theory, which turns out to fit very well with the localization procedure on 
the geometric rough path space with Besov type topology.

In Section 6 we 
prove the large deviation upper estimate.
This part is not very different from the corresponding part of  \cite{in2}.
(However, it is not so easy for those who are not familiar with Watanabe's 
distributional Malliavin calculus).
The key point is the integration by parts formula for Watanabe distributions,
combined with Kusuoka-Stroock's quantitative proof of non-degeneracy of 
Malliavin covariance matrix. 
In Section 7,  using Lyons' continuity theorem and the contraction principles for LDPs,
we prove the LDP conjectured in \cite{tw}
as a simple corollary of our main theorem.

%%\newpage
%%%%%%%%%%%%%%%%%%%%%%%%%%%%%%%%%%%%%%%%%%%%%%%%%%%%%%
%%%%%%%%%%%%%%%%%%%%%%%%%%%%%%%
\section{Setting and Main results}
%%%%%%%%%%%%%%%%%%%%%%%%%%%%%%%%%%%%%%%%%%%%%%%%%%%%%%%%
%%%%%%%%%%%%%%%%%%%%%%%%%%%%%%%

In this section we introduce our setting and state our main results.
Although the setting may seem complicated at first sight, 
we believe that the reader will gradually find it quite natural.
The SDEs we consider in this paper and
our assumptions on the coefficient vector fields are standard.
Our explanation in this section may not be so detailed, but
we will give precise definitions and detailed explanations in later sections.

Let ${\cal W} =C_0 ([0,1], {\mathbb R}^d)$ be the set
of the continuous functions from $[0,1]$ to ${\mathbb R}^d$
which start at $0$. This is equipped with the usual sup-norm.
We denote by ${\cal H}$ and $\mu$
the Cameron-Martin subspace of ${\cal W}$ and the Wiener measure on ${\cal W}$, respectively.
The triple $({\cal W}, {\cal H}, \mu)$ is called the classical Wiener space.
The canonical realization of Brownian motion is denoted by $(w_t)_{0 \le t \le 1}$.

Let $V_{i}: {\mathbb R}^n \to {\mathbb R}^n$ be a vector field 
with sufficient regularity ($0 \le i \le d$).
Precisely, we assume the following regularity condition.
We say that $V_{i}~(0 \le i \le d)$ satisfies Assumption {\bf (A1)}  if
%
%
%(Note that $V_{i}$ itself may have linear growth):
\\
\\
{\bf (A1)}:  $V_{i}$ is of $C^{\infty}$ with bounded derivatives of all order $\ge 1$.
\\
\\
Note that $V_{i}$ itself may have linear growth in the above condition.
When $V_{i}$ is also bounded,  $V_{i}$ is said to be of $C^{\infty}_b$.
($C^{k}_b$ is similarly defined for $k=1,2,\ldots$).

Let $\ve \in (0,1]$ be a small parameter.
Under {\bf (A1)}, we consider
the following SDE of Stratonovich type:
\begin{equation}\label{sc_sde.def}
dX^{\ve}_t = \ve \sum_{i=1}^d  V_i ( X^{\ve}_t) \circ  dw_t^i  + \ve^2  V_0 (X^{\ve}_t)   dt
\qquad
\qquad
\mbox{with \quad $X^{\ve}_0 =x \in {\mathbb R}^n$.}
\end{equation}
When necessary, we will write $X^{\ve}_t = X^{\ve}(t, x, w)$ or $X^{\ve}(t, x)$
and sometimes write $\lambda^{\ve}_t =\ve^2 t$.
Recall that {\bf (A1)} is a standard assumption in Malliavin calculus, 
under which
$X_t^{\ve}$ is ${\bf D}_{\infty}$-functional for all $t \ge 0$ and $\ve \in (0,1]$.

Next we impose a non-degeneracy assumption on the vector fields. 
We set  
$$
\Sigma_1 =\{ V_i ~|~ 1 \le i \le d\} 
\qquad
\mbox{and} 
\qquad
\Sigma_k =\{ [V_i, W] ~|~ 1 \le i \le d, W \in \Sigma_{k-1}\} 
$$
for $k \ge 2$ recursively.
For $x \in {\mathbb R}^n$, we write $\Sigma_k (x) =\{ W(x) ~|~  W \in \Sigma_{k}\}$,
which is a finite subset of ${\mathbb R}^n \cong T_x {\mathbb R}^n$ (i.e., the tangent space at $x$).
We assume the following strong H\"ormander condition everywhere.
\\
\\
{\bf (A2)}:
For any $x \in {\mathbb R}^n$,  $\cup_{k=1}^{\infty}\Sigma_k (x)$ spans 
${\mathbb R}^n \cong T_x {\mathbb R}^n$ in the sense of linear algebra.
\\
\\
Note that the drift vector field $V_0$ is not involved in {\bf (A2)}.

In this paragraph, we will assume $t>0$ and $\ve \in (0,1]$.
It is well-known 
that under {\bf (A1)}--{\bf (A2)},
$X_t^{\ve}$ is non-degenerate in the sense of Malliavin.
Hence, the composition $T (X_t^{\ve}) =T\circ X_t^{\ve}$
is well-defined as a Watanabe distribution (i.e., a generalized Wiener functional) on ${\cal W}$
for any tempered Schwartz distribution $T$ on ${\mathbb R}^n$.
In particular, 
the heat kernel $p^{\ve}_t (x,x')$ (i.e., the density of the law of 
$X^{\ve}(t, x)$ with respect the Lebesgue measure $dx'$)
exists and is equal to ${\mathbb E}[\delta_{x'} ( X^{\ve}(t, x))]$,
where $\delta_{x'} ( X^{\ve}(t, x))$ 
is Watanabe's pullback of the delta function
and ${\mathbb E}$ stands for the generalized expectation.
It is known that 
$p^{\ve}_t (x,x') >0$ for all $x, x' \in {\mathbb R}^n$.
(To check this positivity under our assumptions {\bf (A1)}--{\bf (A2)}, 
combine Theorem 3.41, Aida, Kusuoka, and Stroock \cite{aks} and Theorem 5.3, Kunita \cite{kun} for example.)

Now we introduce the skeleton ODE which corresponds to SDE (\ref{sc_sde.def}).
For a Cameron-Martin path $h \in {\cal H}$,
we consider the following controlled ODE.
\begin{equation}\label{ode.def}
d\phi_t = \sum_{i=1}^d  V_i ( \phi_t) dh_t^i 
\qquad
\qquad
\mbox{with \quad $\phi_0 =x \in {\mathbb R}^n$.}
\end{equation}
Note that this ODE has a unique global solution for any given $h$ under {\bf (A1)}.
The solution will often be denoted by $\phi_t(h),~ \phi(t, x, h)$, etc.
Note the absence of the drift term in (\ref{ode.def}).
%
%We set ${\cal K}^{x, x'} =$ for $x, x' \in {\mathbb R}^n$.

Let ${\cal V}$ be an $l$-dimensional linear subspace of ${\mathbb R}^n$ ($1 \le l \le n$)
and $\Pi_{{\cal V}}: {\mathbb R}^n \to {\cal V}$ be the orthogonal projection.
(For our purpose,  we may and sometimes will 
assume without loss of generality that ${\cal V} = {\mathbb R}^l \times \{ {\bf 0}_{n-l}\}$, 
where ${\bf 0}_{n-l}$ is the zero vector of ${\mathbb R}^{n-l}$.)
Set 
$Y^{\ve}_t = \Pi_{{\cal V}} (X^{\ve}_t)$, which will often be denoted by $Y^{\ve}(t, x, w)$,
and
$\psi(t, x, h) = \Pi_{{\cal V}}\phi(t, x, h)$, 
where $\phi$ is the solution of ODE (\ref{ode.def}).
For $x, x' \in {\mathbb R}^n$ and $a \in {\cal V}$,
define 
${\cal K}^{x, x'} = \{h \in {\cal H} ~|~ \phi (1,x,h) =x' \}$
and 
\[
{\cal M}^{x, a} =
\{ h \in {\cal H} ~|~ \psi(1, x, h) =a\}
=
\bigcup \bigl\{  {\cal K}^{x, x'} ~|~ x' \in \Pi_{{\cal V}}^{-1} (a)  \bigr\}.
\]
By the controllability of ODE (\ref{ode.def}) under {\bf (A2)} (and {\bf (A1)}), 
${\cal K}^{x, x'} \neq \emptyset$ for any $x, x'$.
(See Theorem 5.3, Kunita \cite{kun})
Hence, 
${\cal M}^{x, a} \neq \emptyset$ for any $x, a$.

Let $\lambda (X^{\ve}_t)$ and $\lambda (Y^{\ve}_t)$ be the smallest eigenvalue of 
the Malliavin covariance matrix of $X^{\ve}_t$ and $Y^{\ve}_t$, respectively.
%
%
%Note that since $\Pi_{{\cal V}}^*$ is an isometry, $\lambda (X^{\ve}_t) \le \lambda (Y^{\ve}_t)$.
%
%
It is known that 
$\lambda (X^{\ve}_t)^{-1}$ has moments of all order
(See Nualart \cite{nu} for example.
This is in fact stronger than non-degeneracy of $X^{\ve}_t$ in the sense of Malliavin).
Since $\Pi_{{\cal V}}^*$ is an isometry, $\lambda (X^{\ve}_t) \le \lambda (Y^{\ve}_t)$.
Hence, 
$\lambda (Y^{\ve}_t)^{-1}$ also has moments of all order
and consequently $Y^{\ve}_t$ is non-degenerate in the sense of Malliavin.
(In a similar way,
non-degeneracy of the deterministic Malliavin covariance of 
$\phi(t, x, h)$ at $h \in {\cal H}$
implies that of $\psi(t, x, h)$.)

Therefore, for any $a \in {\cal V}$, $\delta_a (Y^{\ve}_t)$
is a positive Watanabe distribution and 
equal to $(\delta_a \circ \Pi_{{\cal V}}) (X^{\ve}_t)$.
By the positivity of $p_t^{\ve}(x,x')$, 
we can easily see that
${\mathbb E}[ \delta_a (Y^{\ve}_t)] >0$ for all $\ve \in (0,1]$, $t >0$, $x \in {\bf R}^n$, and $a \in {\cal V}$. 
By Sugita's theorem \cite{su}, 
the 
positive Watanabe distribution $\delta_a (Y^{\ve}_1)$ at time $t=1$
is in fact a finite Borel measure on ${\cal W}$,
which will be denoted by $\theta^{\ve}_{x,a}$.

%%%%%%%%%%%%%%%%%%%%%%%%%%%%%%
%%%%%%%%%%%%%%%%%%%%%%%%%%%%%%
%%\vspace{15mm}
%%%%%%%%%%%%%%%%%%%%%%%%%%%%%%
%%%%%%%%%%%%%%%%%%%%%%%%%%%%%%

From here we discuss rough path theory. 
In this paper, we consider the geometric rough path space  
$G\Omega^B_{\alpha, 4m} ( {\mathbb R}^d) $
with Besov-type topology.
We will always assume that the Besov parameter $(\alpha, 4m)$ satisfies the following assumption
so that basic results in \cite{in2} are available;
%(below ${\mathbb N}_{>0} =\{1,2,\ldots\}$);
%
%
\begin{equation}\label{eq.amam}
\frac13 <\al < \frac12, \quad m =1,2,3,\ldots, \quad
\al - \frac{1}{4m} > \frac13, \quad
\mbox{ and } \quad 4m (\frac12 -\al)   >1.
\end{equation}
Under (\ref{eq.amam}),
$G\Omega^B_{\alpha, 4m} ( {\mathbb R}^d) $ is continuously embedded in 
the geometric rough path space  
$G\Omega^H_{\alpha'} ( {\mathbb R}^d) $ with H\"older topology
with $\alpha' = \al - 1/(4m)$.
Intuitively, $\al$ is something like the H\"older exponent 
and $4m$ is a very large even integer.

Let ${\cal L}: {\cal W} \to G\Omega^B_{\alpha, 4m} ( {\mathbb R}^d)$ be
the rough path lift map via the dyadic polygonal approximations, 
which is defined outside a slim subset of ${\cal W}$
and $\infty$-quasi continuous.
(When this is regarded as a geometric rough path space-valued random variable,
we will often write ${\bf W} = {\cal L}(w)$. 
On the other hand,
a generic element of 
the geometric rough path space is denoted by ${\bf w}$, which is not random.)
Since 
${\cal L}$ is quasi-surely defined, we can lift 
the measure $\theta^{\ve}_{x,a}$ to a
measure on $G\Omega^B_{\alpha, 4m} ( {\mathbb R}^d)$.
We write 
$\mu^{\ve}_{x,a} =(\ve \cdot {\cal L})_* [\theta^{\ve}_{x,a}]$,
where the "dot" 
stands for the dilation on $G\Omega^B_{\alpha, 4m} ( {\mathbb R}^d)$.
We denote by $\hat\theta^{\ve}_{x,a}$ and $\hat\mu^{\ve}_{x,a}$
the normalized measure of 
$\theta^{\ve}_{x,a}$ and $\mu^{\ve}_{x,a}$, respectively.
(Since the total mass of $\theta^{\ve}_{x,a}$ or of $\mu^{\ve}_{x,a}$ 
equals ${\mathbb E}[ \delta_a (Y^{\ve}_t)] >0$, this normalization is well-defined.)

%%%%%%%%%%%%%%%%%%%%%%%%%%%%%%
%%%%%%%%%%%%%%%%%%%%%%%%%%%%%%
%%\vspace{15mm}
%%%%%%%%%%%%%%%%%%%%%%%%%%%%%%
%%%%%%%%%%%%%%%%%%%%%%%%%%%%%%

Set a rate function $I_1 : G\Omega^B_{\al, 4m} ({\bf R}^n) \to [0, \infty]$ as follows;
\begin{align}
I_1 ({\bf w}) 
= 
\begin{cases}
    \|h\|^2_{{\cal H}}/2 & (\mbox{if ${\bf w}= {\cal L}(h)$ for some $h \in {\cal M}^{x,a}$}), \\
    \infty &  (\mbox{otherwise}).
  \end{cases}
\nn
\end{align}
This rate function $I_1$ is actually good.
We also set $\hat{I}_1  ({\bf w}) 
= I_1 ({\bf w}) - \min\{ \|h\|^2_{{\cal H}}/2  ~|~h \in {\cal M}^{x,a}\}$.
Note that the minimum above exists.

The following theorem is our main result in this paper.
It states that the family of finite measures
$\{\mu^{\ve}_{x,a} \}_{0 <\ve \le 1}$
satisfies 
an LDP of Schilder type on $G\Omega^B_{\alpha, 4m} ( {\mathbb R}^d)$
as $\ve \searrow 0$.

\begin{tm}\label{tm.main}
Assume {\bf (A1)} and {\bf (A2)} and the condition (\ref{eq.amam}).
Then, we have the following {\rm (i)}--{\rm (ii)}:
\\
\noindent
{\rm (i)}~ The family $\{ \mu^{\ve}_{x,a}\}_{\ve >0}$ 
of finite measures
satisfies an LDP 
on $G\Omega^B_{\al, 4m} ({\bf R}^d)$ as $\ve \searrow 0$ with a good rate function $I_1$, that is, 
for any Borel set $A \subset G\Omega^B_{\al, 4m} ({\bf R}^d)$,
the following inequalities hold;
\begin{align}
- \inf_{{\bf w} \in A^{\circ} } I_1 ({\bf w})  
\le
 \liminf_{\ve \searrow 0 } \ve^2 \log \mu^{\ve}_{x,a} (A)
\le
 \limsup_{\ve \searrow 0 } \ve^2 \log \mu^{\ve}_{x,a} (A)
\le 
- \inf_{{\bf w} \in \bar{A} } I_1  ({\bf w}).
\nn
\end{align}
\noindent
{\rm (ii)}~ 
The family $\{ \hat\mu^{\ve}_{a,a'}\}_{\ve >0}$ 
of probability measures 
satisfies an LDP 
on $G\Omega^B_{\al, 4m} ({\bf R}^d)$ as $\ve \searrow 0$ with a good rate function $\hat{I}_1$.
\end{tm}

Since the whole set is both open and closed,
Theorem \ref{tm.main}, {\rm (i)} implies that 
$$
\lim_{\ve \searrow 0 } \ve^2 \log \mu^{\ve}_{x,a} ( G\Omega^B_{\al, 4m} ({\bf R}^d) ) 
= \lim_{\ve \searrow 0 } \ve^2 \log {\mathbb E}[ \delta_a (Y^{\ve}_t)] 
=
- \min\{ \|h\|^2_{{\cal H}}/2  ~|~h \in {\cal M}^{x,a}\}.$$
Therefore, 
Theorem \ref{tm.main}, {\rm (ii)} is immediate from {\rm (i)}.

Theorem \ref{tm.main} above also holds with respect to 
$\alpha'$-H\"older geometric rough path topology for any $\al' \in (1/3, 1/2)$,
because we can find $\al, m$ with (\ref{eq.amam}) 
such that $(\al, 4m)$-Besov topology is stronger than $\alpha'$-H\"older topology.

%%%%%%%%%%%%%%%%%%%%%%%%%%%%%%
%%%%%%%%%%%%%%%%%%%%%%%%%%%%%%
%%\vspace{15mm}
%%%%%%%%%%%%%%%%%%%%%%%%%%%%%%
%%%%%%%%%%%%%%%%%%%%%%%%%%%%%%

From the contraction principle for LDPs, 
it is obvious that, for any continuous map $F$ from the geometric rough path space 
to a Hausdorff topological space, 
the image measure $F_* [\hat\mu^{\ve}_{x,a}]$ satisfies an LDP, too.
As an example of such continuous maps, 
we may take a Lyons-It\^o map associated with coefficient vector fields 
which are different from $V_i$'s.

Let $A_{i}: {\mathbb R}^N \to {\mathbb R}^N$ be a vector field 
which satisfies {\bf (A1)} ($0 \le i \le d$).
Note that they may be different from $V_i$.
For $\ve \in (0,1]$, 
we also consider
the following SDE of Stratonovich type; 
\begin{equation}\label{sc_sde_A.def}
dZ^{\ve}_t = \ve \sum_{i=1}^d  A_i ( Z^{\ve}_t) \circ  dw_t^i  + \ve^2  A_0 (Z^{\ve}_t)   dt
\qquad
\qquad
\mbox{with \quad $Z^{\ve}_0 =z \in {\mathbb R}^N$.}
\end{equation}
For $h \in {\cal H}$,
we consider the following controlled ODE;
\begin{equation}\label{ode_z.def}
d\zeta_t = \sum_{i=1}^d  A_i ( \zeta_t) dh_t^i 
\qquad
\qquad
\mbox{with \quad $\zeta_0 =z \in {\mathbb R}^n$.}
\end{equation}
We may write $Z^{\ve}_t = Z^{\ve}(t, z, w)$ or $\zeta_t =\zeta_t(h) = \zeta (t, z, h)$, etc.
We denote by $\tilde{Z}^{\ve}=\tilde{Z}^{\ve}(\,\cdot\,, z, w)$
an $\infty$-quasi continuous modification of 
\begin{equation}\label{tildeZ.def}
{\cal W} \ni w \mapsto Z^{\ve}(\,\cdot\,, z, w) \in C^{\alpha -H}([0,1], {\mathbb R}^N)
\qquad\qquad
(1/3 <\alpha<1/2).
\end{equation}
Here, the set on the right hand side stands for the space of $\alpha$-H\"older 
continuous paths in ${\mathbb R}^N$.
Since $\tilde{Z}^{\ve}$ is defined uniquely up to a slim subset of ${\cal W}$,
the pushforward measures of $\theta^{\ve}_{x,a}$ 
and $\hat\theta^{\ve}_{x,a}$  by the map $\tilde{Z}^{\ve}$ are well-defined.

As a corollary of our main theorem, we can prove
an LDP as $\ve \searrow 0$ for these measures. 
Before stating it, let us first define good rate functions 
$I_2, \hat{I}_2 : C^{\alpha -H}([0,1], {\mathbb R}^N) \to [0,\infty]$.
Set
\begin{align}
I_2 (b) 
= 
\begin{cases}
  \inf\{  \|h\|^2_{{\cal H}}/2 ~|~ \mbox{
  $h \in {\cal M}^{x,a}$ such that $b = \zeta(\,\cdot\,,z,h)$ } \},
    \\
    \infty, \qquad \mbox{(if no $h \in {\cal M}^{x,a}$ satisfies that $b = \zeta(\,\cdot\,,z,h)$)}.
  \end{cases}
\nn
\end{align}
and $\hat{I}_2 (b) 
= I_2 (b) - \min\{ \|h\|^2_{{\cal H}}/2  ~|~h \in {\cal M}^{x,a}\}$.

%%%%%%%%%%%%%%%%%%%%%%%%%%%%%%
%%%%%%%%%%%%%%%%%%%%%%%%%%%%%%
%%\vspace{15mm}
%%%%%%%%%%%%%%%%%%%%%%%%%%%%%%
%%%%%%%%%%%%%%%%%%%%%%%%%%%%%%

\begin{co}\label{co.TW}
Let $1/3 <\al <1/2$. 
Assume {\bf (A1)} for both $V_i$ and $A_i~(0 \le i \le d)$ 
and assume {\bf (A2)} for $V_i~(0 \le i \le d)$.
Then, we have the following {\rm (i)}--{\rm (ii)}:
\\
\noindent
{\rm (i)}~ The family $\{ \tilde{Z}^{\ve}(\,\cdot\,, z)_* [\theta^{\ve}_{x,a}]\}_{\ve >0}$ 
satisfies an LDP
on $C^{\alpha -H}([0,1], {\mathbb R}^N)$ as $\ve \searrow 0$ with a good rate function $I_2$, that is, 
for any Borel set $A \subset C^{\alpha -H}([0,1], {\mathbb R}^N)$,
the following inequalities hold;
\begin{align}
- \inf_{ b \in A^{\circ} } I_2 ( b )  
&\le
 \liminf_{\ve \searrow 0 } \ve^2 \log 
 \theta^{\ve}_{x,a} (\{w \in {\cal W}~|~ \tilde{Z}^{\ve}(\,\cdot\,, z,w) \in A\})
\nn\\
&\le
 \limsup_{\ve \searrow 0 } \ve^2 \log \theta^{\ve}_{x,a}
  (\{w \in {\cal W}~|~ \tilde{Z}^{\ve}(\,\cdot\,, z,w) \in A\})
\le 
- \inf_{ b \in \bar{A} } I_2 ( b).
\nn
\end{align}
\noindent
{\rm (ii)}~The family $\{ \tilde{Z}^{\ve}(\,\cdot\,, z)_* [\hat\theta^{\ve}_{x,a}]\}_{\ve >0}$ 
of probability measures
satisfies an LDP
on $C^{\alpha -H}([0,1], {\mathbb R}^N)$ as $\ve \searrow 0$ with a good rate function $\hat{I}_2$.
\end{co}

\begin{re}\label{re.HB_TW}
In the formulation of Corollary \ref{co.TW} above, 
H\"older path space $C^{\alpha -H}([0,1], {\mathbb R}^N)$ is used,
while Besov-type path space is used in Theorem 2.1, p. 200, \cite{tw}. 
However, by adjusting H\"older/Besov parameters, 
we see that these two formulations are equivalent.
Therefore, 
Corollary \ref{co.TW} above is equivalent to the LDP
conjectured in Theorem 2.1, \cite{tw}.
\end{re}

Corollary \ref{co.TW} above immediately implies an LDP of Freidlin-Wentzell type 
for pinned diffusion measures as follows.

Take $n=l =N$, $x=z$ and $V_i =A_i$ for all $i$. 
We write $a=x' \in {\mathbb R}^n$.
Then, $X^{\ve}_t = Y^{\ve}_t= Z^{\ve}_t$, $\phi_t = \psi_t =\zeta_t$, and ${\cal M}^{x,a}={\cal K}^{x,x'}$.
In this case, 
$ \tilde{Z}^{\ve}(\,\cdot\,, z)_* [\hat\theta^{\ve}_{x,a}]$ is nothing but 
the pinned diffusion measure $Q^{\ve}_{x,x'}$ 
associated to the generator $\ve^2 \{V_0 +(1/2)\sum_{i=1}^d V_i^2 \}$
(or equivalently, to the heat kernel $p^{\ve}_t$)
with the starting point $x$ and the ending point $x'$.

Then, we have the following result. 
The proof is almost obvious.
\begin{co}\label{co.pFW}
Let $1/3 <\al <1/2$ and assume {\bf (A1)} and {\bf (A2)}.
The family $\{Q^{\ve}_{x,x'} \}_{\ve >0}$ 
satisfies an LDP 
on $C^{\alpha -H}([0,1], {\mathbb R}^N)$ as $\ve \searrow 0$ with a good rate function $I_2^{\prime}$.
Here, $I_2^{\prime}$ is given by 
\begin{align}
 I_2^{\prime} (b) 
= 
\begin{cases}
  \inf\{  \|h\|^2_{{\cal H}}/2 ~|~ \mbox{
  $h \in {\cal K}^{x,x'}$ such that $b = \phi(\,\cdot\,,z,h)$ } \}  - \min\{ \|h\|^2_{{\cal H}}/2  ~|~
  h \in {\cal K}^{x,x'}\},
    \\
    \infty, \qquad \mbox{(if no $h \in {\cal K}^{x,x'}$ satisfies that $b = \phi(\,\cdot\,,z,h)$)}.
  \end{cases}
\nn
\end{align}
\end{co}

We remark that Bailleul proved 
an LDP parallel to Corollary \ref{co.pFW} on  compact manifolds
in \cite{bai} (and in its extended version \cite{bmn} with Mesnager and Norris).
Their method is basically  analytic (with a little bit of rough path theory) and different from ours.
Their result can be viewed as a hypoelliptic version of Hsu's result in \cite{hsu}
for pinned Brownian motions on compact Riemannian manifolds.

%%%%%%%%%%%%%%%%%%%%%%%%%%%%%%
%%%%%%%%%%%%%%%%%%%%%%%%%%%%%%
%%

\begin{re}\label{re.weak}
One cannot replace the "strong H\"ormander" condition in Theorem 2.1 
by the "H\"ormander" condition. 
We have the following counterexample.
Consider the following two-dimensional SDE driven by one-dimensional 
Brownian motion.
\[
dX^{\ve, 1}_t = \ve dw_t,  \qquad  dX^{\ve, 2}_t = \ve^2 X^{ \ve, 1}_t dt.
\]
The coefficient vector fields satisfy the H\"ormander condition everywhere, 
but nowhere the strong  H\"ormander condition.
If the solution starts at the origin, the law of $(X_{1}^{\ve,1},X_{1}^{\ve,2})$ 
is the centered Gaussian measure 
with the covariance
\[
\begin{pmatrix}
\ve^2 & \ve^4/2 \\
 \ve^4/2 & \ve^6/3
 \end{pmatrix}.
\]
Then, it is easy to see that
$p^{\ve}_1 ((0,0), (0,x^2)) = \sqrt{3} (\pi \ve^4)^{-1} \exp ( -(x^2)^2 /(6 \ve^6))$.
If $x^2 \neq 0$, then 
$\lim_{\ve \searrow 0} \ve^2 \log p^{\ve}_1 ((0,0), (0,x^2)) = -\infty$.
On the other hand, we have
 ${\cal K}^{(0,0), (0,x^2)} \neq \emptyset$.
Therefore, the heat kernel does not behave in the way described in Theorem \ref{tm.main}.
(Recall that the heat kernel is the weight of the whole set in our setting).
 \end{re}

\begin{re}\label{re.prev}
Loosely speaking, our main results above  generalize the ones for
the elliptic case in the author's previous paper \cite{in2}.
However, the results in this paper do not cover all of the results in \cite{in2}
for the following reasons:
\\
{\rm (i)}~In this paper the strong  H\"ormander condition is assumed at any point,
while in \cite{in2}, the ellipticity condition is assumed only at the starting point
and at some point vector fields may even be degenerate (i.e., do not even satisfy the H\"ormander condition).
\\
{\rm (ii)}~
In \cite{in2} the drift vector field is of the form $V_0 (\ve, x)$ 
and is quite general.
However, it is of the form $\ve^2 V_0 (x)$ in this paper.
Although it may be possible to generalize our results 
for a drift term of the form $V_0 (\ve, x)$ with $V_0 (0, x) \equiv 0$,
it is probably impossible 
if $V_0 (0, x)$ do not vanish identically. 
(This guess is based on an observation of small noise asymptotics of the heat kernel
in Section 3, Ben Arous and L\'eandre \cite{bal}).
\end{re}

%%%%%%%%%%%%%%%%%%%%%%%%%%%%%%%%%%%%%%%%%%%%%%%%%%%%%%%%
%%%%%%%%%%%%%%%%%%%%%%%%%%%%%%%\newpage
\section{Skeleton ODE}

In this section we study the solution $\phi_t (h)= \phi (t, x, h)$ of 
the skeleton ODE (\ref{ode.def}).
Note that it always has a global solution under {\bf (A1)}.
The aim of this section 
is to prove that a Fr\'echet differentiable map
$h \mapsto \phi (1,x,h)$ is non-degenerate at sufficiently many $h$'s 
under strong H\"ormander condition on the vector fields. (See Proposition \ref{pr.nondeg.dmc}).
It will play a crucial role in the lower estimate for the LDP
in our main theorem (Theorem \ref{tm.main}).
We emphasize again that 
the absence of the drift term in (\ref{ode.def}) has a significant meaning
and many parts of this section would fail if (\ref{ode.def}) had a drift term.

\subsection{Basic properties of skeleton ODE}

First we set some notations.
For $T>0$, ${\cal H}_T$ denotes ${\mathbb R}^d$-valued 
Cameron-Martin space on the time interval $[0,T]$, that is, 
\[
{\cal H}_T = \bigl\{ h: [0,T]\to {\mathbb R}^d ~|~ \mbox{$h =\int_0^{\cdot} \dot{h}_s ds$
for some $\dot{h} \in L^2([0,T], {\mathbb R}^d)$} \bigr\}.
\]
The Hilbert norm is naturally defined by 
$\| h\|_{{\cal H}_T } =  \| \dot{h}\|_{L^2[0,T]}$ as usual.
When $T=1$, we simply write ${\cal H} ={\cal H}_1$.

For $h \in {\cal H}_T$, the reversed path 
$\overline{h} \in  {\cal H}_T$ is defined by $\overline{h}_t =h_{T -t} -h_T$.
Concatenation of $h \in {\cal H}_T$ and $k \in {\cal H}_S$
is denoted by $h*k \in {\cal H}_{T+S}$,
which is defined by $(h*k) (t) =h(t)$ for $0 \le t \le T$
and 
$(h*k) (t) =k_{t-T} +h_T$ for $T \le t \le T+S$

%For $x, x' \in {\mathbb R}^n$, we set ${\cal K}^{x, x'}=\{ h \in {\cal H} ~|~ \phi(1, x, h) =x'\}$.
%
Thanks to {\bf (A2)}, ODE (\ref{ode.def}) is strongly completely controllable
(Theorem 5.3, Kunita \cite{kun}).
Hence, for any $x, x'$ and $T >0$, 
there exists $h \in {\cal H}_T$ 
such that $\phi(T,x, h) =x'$ and, in particular, 
${\cal K}^{x, x'} := \{ h \in {\cal H} ~|~ \phi(1, x, h) =x'\}\neq \emptyset$.

%%%%%%%%%%%%%%%%%%%%%%%%%%%%%%%%%%%%%%%%%%%%%%%%%%%%%%%%
%%  deterministic  Jacobian process
%%%%%%%%%%%%%%%%%%%%%%%%%%%%%%%%%%%%%%%%%%%%%%%%%%%%%%%%%

Now we introduce Jacobian ODE of (\ref{ode.def}) and its inverse.

\begin{align}
dJ_t &= \sum_{i=1}^d  \nabla V_i ( \phi_t) J_t dh_t^i 
&\mbox{with \quad $J_0 = {\rm Id}_n$,}
\label{ode_J.def}
\\
dK_t &= - \sum_{i=1}^d  K_t \nabla V_i ( \phi_t)  dh_t^i 
&\mbox{with \quad $K_0 = {\rm Id}_n$.}
\label{ode_K.def}
\end{align}
Here, $J, K, \nabla V_i $ are all $n \times n$ matrices.
Note that $K_t =J_t^{-1}$.
When dependency on $h$ and $x$ needs to be specified, 
we write $J_t (h)$ or $J(t,x,h)$, etc.

The map $h \in {\cal H} \mapsto \phi_t(h) = \phi (t,x,h)  \in {\mathbb R}^n$
is of Fr\'echet-$C^1$ for each $t \in [0,1]$ and $x \in {\mathbb R}^n$.
The Fr\'echet derivative $D \phi_t(h) \in L({\cal H}, {\mathbb R}^n)$
is explicitly given by
\begin{equation}\label{Dphi.eq}
D \phi_t(h) \la  k \ra 
=
J_t (h) \sum_{i=1}^n  \int_0^t  J_s (h)^{-1} V_i ( \phi_s (h)) \dot{k}^{i}_s ds
\end{equation}
The deterministic Malliavin covariance of $\phi_1$ at $h$
(and at time $t=1$)
 is defined by 
\begin{equation}\label{dtmmalcov.def}
\sigma_{\phi_1} (h) = \Bigl( \la D \phi_1^i (h) ,  D \phi_1^j (h)  \ra_{{\cal H}^*}  \Bigr)_{1 \le i,j \le n}
=
D \phi_1 (h) \circ [ D \phi_1 (h) ]^*,
\end{equation}
where the superscript $*$ stands for the adjoint operation. 
From (\ref{Dphi.eq}) and (\ref{dtmmalcov.def}) we can easily see that 
$\sigma_{\phi_1} (h) = J_1(h)   C(h) J_1(h)^*$ with
\begin{equation}\label{rep.dmc.eq}
C(h) 
=
\int_0^1  J_s(h)^{-1}  {\bf V} (\phi_s(h))   {\bf V} (\phi_s(h))^* J_s(h)^{-1, *}ds
\end{equation}
Here, we set $ {\bf V} (x) =[V_1(x), \ldots, V_d(x) ]$, which is an $n \times d$ matrix.
Note that the surjectivity of  the linear map $D \phi_1 (h):{\cal H} \to {\mathbb R}^n$ is 
equivalent to non-degeneracy of 
the deterministic Malliavin covariance $\sigma_{\phi_1} (h) $, 
which in turn is equivalent to 
non-degeneracy of  $C(h)$ since $J_1(h)$ is always invertible.

The following is the main result in this section.
Unlike in the elliptic case, there exists $h$ such that $\sigma_{\phi_1} (h)$ is degenerate.
(For example, think of the constant path $0 \in {\cal H}$.)
However, there are sufficiently many $h$'s for which $\sigma_{\phi_1} (h)$ is non-degenerate.
The precise statement is given as follows.
\begin{pr}\label{pr.nondeg.dmc}
Assume {\bf (A1)} and {\bf (A2)}. 
Let $x, x' \in {\mathbb R}^n$ and $h \in {\cal K}^{x, x'}$ be arbitrary. 
Then, we have the following;
\\
\noindent
{\rm (i)}~For any $\ve >0$, there exists $h^{\ve} \in {\cal K}^{x, x'}$ 
such that $\|  h -h^{\ve} \|_{{\cal H}} <\ve$ and
$\sigma_{\phi_1} (h^{\ve})$ is non-degenerate.
\\
\noindent
{\rm (ii)}~Moreover, $h^{\ve}$ in {\rm (i)} above can be chosen so that 
$\la h^{\ve}, \,\cdot\, \ra_{{\cal H}}$  naturally extends to 
a continuous linear functional on the  Wiener space ${\cal W}$.
\end{pr}

The proof of Proposition \ref{pr.nondeg.dmc}, {\rm (i)} 
will be given in the subsequent subsections.
Once we have Proposition \ref{pr.nondeg.dmc}, {\rm (i)},  
we can prove {\rm (ii)} by using the following lemma (with 
${\cal K} ={\cal H}, ~ {\cal L}={\cal W}^*$ and
${\cal W}^*\subset {\cal H}^* \cong {\cal H}$).
%
%

%%%%%%%%%%%%%%%%%%%%%%%%%%%%%%%

\begin{lm}\label{lm.inver}
Let ${\cal K}$ be a real Hilbert space and $\xi \in {\cal K}$.
Assume that {\rm (i)}~ $F$ is an ${\mathbb R}^n$-valued Fr\'echet-$C^1$ map 
 defined on a neighborhood of $\xi$ with a bounded derivative $D F$
and {\rm (ii)}~$DF (\xi): {\cal K} \to {\mathbb R}^n$ is a surjective linear map.
Let ${\cal L}$ be a real Banach space which is continuously and densely 
embedded in ${\cal K}$.
Then, there exists $\xi_j \in {\cal L}~(j=1,2,\ldots)$ such that 
$\lim_{j\to \infty} \| \xi_j- \xi \|_{\cal K} =0$ and $F(\xi_j) = F(\xi)$ for all $j$.
(Necessarily, $DF (\xi_j)$ is also surjective for large enough $j$.)
\end{lm}

\Proof
This lemma was proved in \cite{in2}.
\QED

%%%%%%%%%%%%%%%%%%%%%%%%%%%%%%%

Before closing this subsection, 
we prove two simple lemmas for later use. 
For $h \in {\cal H}$ and $x \in {\mathbb R}^n$ and a vector field $W: {\mathbb R}^n \to {\mathbb R}^n$, 
we set 
$Q^W_t = J_t^{-1} W (\phi_t)$.
Note that both $J^{-1}$ and $\phi$ depend on $h$ and $x$.
(We will sometimes write $Q^W_t(h)$ or $Q^W (t, x, h)$, etc.)

\begin{lm}\label{lm.abc}
{\rm (i)}~Let $W: {\mathbb R}^n \to {\mathbb R}^n$ be a smooth vector field.
Then, we have
\[
d Q^W_t(h) = \sum_{i=1}^d  Q^{[V_i, W]}_t(h) \dot{h}_t^i dt.
\]
\\
{\rm (ii)}~
For any $v \in {\mathbb R}^n$,  we have 
\[
v^* C(h) v = \sum_{i=1}^d  \int_0^1 |\la v, Q^{V_i}_s (h) \ra|^2 ds.
\]
In particular, 
if for any $v$ with $\|v\|=1$, there exist $t \in [0,1]$ and $i ~(1 \le i \le d)$
such that $\la v, Q^{V_i}_t (h) \ra \neq 0$,
then $C(h)$ and $\sigma_{\phi_1} (h)$ are non-degenerate.
\end{lm}

\Proof
The first assertion can easily be seen from (\ref{ode.def})--(\ref{ode_K.def}). 
The second one is shown by simple calculation of block matrices 
and is a routine.
So, the proof is omitted.
\QED

The next lemma is quite simple. (So we omit a proof.)
However,  note that the absence of a drift term in (\ref{ode.def})--(\ref{ode_K.def})  
is crucially important here.
If they had a drift term, this lemma would fail.
%
%%%%%%  Lemma %%%%
%
\begin{lm}\label{lm.defg}
{\rm (i)}~If $(\phi_t, J_t, K_t)_{0 \le t \le T}$ is the solution of  
ODEs (\ref{ode.def}), (\ref{ode_J.def}), and (\ref{ode_K.def})  driven by $h \in {\cal H}_T$
 with the initial condition $(x, {\rm Id}_n, {\rm Id}_n)$,
then $(\phi_{T-t}, J_{T-t}, K_{T-t})_{0 \le t \le T}$ is the solution of 
(\ref{ode.def}),  (\ref{ode_J.def}), and (\ref{ode_K.def})  driven by 
the reversed path $\bar{h} \in {\cal H}_T$,
with the initial condition $(\phi_T, J_T, K_T)$
\\
{\rm (ii)}~For $h \in {\cal H}_T$,  let $(\phi_t, J_t, K_t)_{0 \le t \le 2T}$
be the solution of  
(\ref{ode.def}), (\ref{ode_J.def}), and (\ref{ode_K.def})  driven by $h * \bar{h} \in {\cal H}_{2T}$,
then $(\phi_{2T}, J_{2T}, K_{2T})=(\phi_0, J_0, K_0)$.
\\
{\rm (iii)}~Let $T >0$ and $\beta \in (0,T)$.
If $(\phi_t, J_t, K_t)_{0 \le t \le T}$ is the solution of  
(\ref{ode.def}), (\ref{ode_J.def}), and  (\ref{ode_K.def})  driven by $h \in {\cal H}_T$
 with the initial condition $(x, {\rm Id}_n, {\rm Id}_n)$,
then 
$$(\phi_{T(t -\beta) /(T -\beta)}, J_{T(t -\beta) /(T -\beta) }, 
K_{ T(t -\beta) /(T -\beta)})_{\beta \le t \le T}
$$
is the solution of  
(\ref{ode.def}), (\ref{ode_J.def}), and (\ref{ode_K.def})  driven by $k$
 with the initial condition $(\phi_{\beta}, J_{\beta}, K_{\beta}) = (x, {\rm Id}_n, {\rm Id}_n)$,
where $k$ is defined by $k_t = h_{ T(t -\beta) /(T -\beta)}$ 
on the time interval  $[\beta, T]$.
\end{lm}

%%%%%%%%%%%   Lemma %%%%%%%%%%%
%
%
\begin{lm}\label{lm.hij}
Let $h \in {\cal H}$. 
Assume that, for any $\beta \in (0,1)$,  
$h^{\beta} \in {\cal H}$ which satisfies  the following property is given:
\[
\mbox{
$|\dot{h}^{\beta}_t | \le 1$ for a.a. $t \in [0, \beta]$ 
\quad
and
\quad
$h^{\beta}_t  = h_{(t -\beta)/(1 -\beta)}$ on $[\beta,1]$.
}
\]
Then,  $h^{\beta} \to h$ in ${\cal H}$ as $\beta \searrow 0$.
\end{lm}

\Proof
Without loss of generality, we may assume $d=1$.
It is sufficient to show that $\dot{h}^{\beta} \to \dot{h}$ in $L^2$-norm.

For any $\ve >0$,  there exists a continuous function $f :[0,1] \to {\mathbb R}$
such that $\|  f - \dot{h} \|_{L^2} <\ve$.
It is easy to see that
\begin{eqnarray*}
\| \dot{h}  - \dot{h}^{\beta} \|_{L^2}  
&\le&
\| \dot{h}  \|_{L^2[0,\beta]}  + \| \dot{h}^{\beta}  \|_{L^2[0,\beta]} 
+
\| \dot{h}  - \dot{h}^{\beta} \|_{L^2[\beta, 1]}.
\end{eqnarray*}
The first and the second terms clearly vanish as $\beta \searrow 0$.
The third term is dominated by 
\[
\| \dot{h} -f   \|_{L^2[\beta,1]}  
+
\Bigl\|  f  -\frac{1}{1 -\beta} f  \Bigl( \frac{\cdot -\beta}{1 -\beta} \Bigr)   \Bigr\|_{L^2[\beta,1]}  
+
\frac{1}{1 -\beta}
\Bigl\|    f \Bigl( \frac{\cdot -\beta}{1 -\beta} \Bigr)   
 -  \dot{h} \Bigl( \frac{\cdot -\beta}{1 -\beta}  \Bigr) \Bigr\|_{L^2[\beta,1]}. 
\]
From the way $f$ is chosen, the sum of  the first and the third term is dominated by 
$3\ve$ if  $\beta>0$ is sufficiently small.
Due to the uniform continuity of $f$,  
the second term vanishes as $\beta \searrow 0$.
Thus, we have shown 
$\limsup_{\beta \searrow 0 }\| \dot{h}  - \dot{h}^{\beta} \|_{L^2}  \le 3\ve$.
Letting $\ve \searrow 0$, we finish the proof of the lemma.
\QED

%%%%%%%%%%%%%%%%%%%%%%%

We will fix an arbitrary initial point $x$ in what follows. 
We say $\{V_i\}$ satisfies strong H\"ormander condition of degree $N$ at $x$ 
if $N$ is the smallest integer such that 
$\cup_{k=1}^N \Sigma_k (x)$ linearly spans ${\mathbb R}^n$. 
In this case there exists a subset $\Lambda \subset \cup_{k=1}^N \Sigma_k$
of cardinality $n$
such that $\{ W (x') ~|~  W \in  \Lambda\}$ 
linearly spans ${\mathbb R}^n$ for any $x'$ which is sufficiently close to $x$.  
(We will write $\Lambda =\{ W_1, \ldots, W_n \}$.)
By compactness,  it holds that
$3\lambda := \inf_{v \in {\mathbb S}^{n-1}} \max_{1 \le j \le n} 
|\la v, W_j (x)  \ra|  >0$,
where
${\mathbb S}^{n-1} = \{ v \in {\mathbb R}^n ~|~  |v| =1 \}$ is the unit sphere.

%%%%

%%%%%%%%%%%\vspace{16mm}

%%%%%%%%%%%%%  Lemma 3.6 %%%%%
%%%%%
\begin{lm}\label{lm.klm}
Keep the same notations as above.
For sufficiently small $T >0$, the following property holds:
For any $v \in {\mathbb S}^{n-1}$, there exists $W \in \Lambda$ such that
\[
\inf
\bigl\{
|\la v, Q_t^{W} (h)  \ra|  
\, \mid \,
 0 \le t \le T, \quad  \mbox{$h \in {\cal H}$ with $|\dot{h}_s| \le 1$ for a.a. $s \in [0,T]$}
\bigr\}
\ge \lambda >0.
\]
\end{lm}

\Proof
Set $E_j =\{ v \in {\mathbb S}^{n-1} \mid  |\la v, W_j (x)  \ra|  \ge 3\lambda \}$
for $1 \le j \le n$.
Then, each $E_j $ is compact and ${\mathbb S}^{n-1}  =\cup_j E_j$.
Since $W_j$ is continuous in $x$,  
there exists $r>0$ such that
\[
\inf
\bigl\{
|\la v, W_j (x')  \ra|
\, \mid \,
x' \in B_r (x), \,  v \in E_j
\bigr\}
\ge 2\lambda >0
\qquad
(1 \le j \le n),
\]
where $B_r (x) =\{ x' \in  {\mathbb R}^n ~|~  |x' -x | <r   \}$ 
is the ball of radius $r>0$ centered at $x$.

Let $v \in E_j$.
If $T$ is sufficiently small, then $\phi_t (h)$ stays inside $B_r (x)$.
Therefore, 
$|\la v, W_j ( \phi_t (h) )  \ra| \ge 2\lambda$ when $0 \le t \le T$.
On the other hand,  there exists a constant
$M>0$ (independent of such an $h$) which satisfies that
$|J_t(h)^{-1} -J_s(h)^{-1} | \le M |t-s|$ for all $s,t \in [0,T]$.
Hence, by taking 
$T \le \lambda M^{-1}  ( \max_j \sup \{ |W_j (x')|   \mid x' \in  B_r (x) \}  )^{-1}$,
we can prove the lemma
since $Q_t^{W_j} (h) = J_t(h)^{-1} W_j ( \phi_t (h) )$.
Notice that the choice of $T$ is independent of $v$. 
\QED

\begin{re}
In what follows, the constants $r, T \in (0,1)$ which appear 
in (the proof of) Lemma \ref{lm.klm} above will be fixed.
(Of course, so will $\lambda >0$.)
\end{re}

%%%%%%%%%%%%%%%%%%%%%%%%%%%%%%%%%%%%%%%%%%%%%%%%%%%%%%%%%%%%%%%%%%
%%%%%%%%%%%%%%%%%%%%%%%%%%%%%%%%%%%%%%%%%%%%%%%%%%%%%%%%%%%%%%%%%%
%%%%%%%%%%%%%%%%%%%%%%%%%%%%%%%%%%%%%%%%%%%%%%%%%%%%%%%%%%%%%%%%%%%%%%%%
%%%%%%%%%%%%%%%%%%%%%%%%%%%%%%%%%%%%%%%%%%%%%%%%%%%%%%%%%%%%%%%%%%

\subsection{Proof of Proposition \ref{pr.nondeg.dmc} {\rm (i)}: For degrees $1, 2, 3$}
In this subsection we prove 
Proposition \ref{pr.nondeg.dmc} {\rm (i)} when $N$, i.e., the degree of hypoellipticity at the 
initial point $x$, is $1, 2$ or $3$.
(Strictly speaking, this subsection is not necessary. 
However, we believe it helps the reader understand what is going on in 
the proof for the general case in the next subsection.)

%%%%

Before doing so,  we set a few notations for general $N$.
First, Let $L>0$ be the smallest constant such that 
\[
| W(x_1)- W(x_2) |
\le L |x_1 -x_2|  
\qquad
\qquad
(x_1, x_2 \in B_r (x),  \,  W \in \cup_{k=1}^N \Sigma_k).
\]
Note that $r$ has already been determined.
From this we can see the following estimate (\ref{M.ineq}):
Let $T>0$ be as in Lemma \ref{lm.klm}. 
Then, there exists $M>0$ such that
for any $s,t \in [0,T]$,  $W \in \cup_{k=1}^N \Sigma_k$ and 
$h$ with $|\dot{h}_t| \le 1$ for a.a. $t \in [0,T]$,  it holds that
\begin{equation}
\label{M.ineq}
| Q_t^{W} (h)- Q_s^{W} (h) | \le M|t-s|.
\end{equation}

For $\tau >0$,  $i \in \{1, \ldots, d\}$, and $\kappa \in \{\pm 1 \}$,
we set $\xi^{\tau, i, \kappa} \in {\cal H}_{\tau}$ by 
$\dot{\xi}^{\tau, i, \kappa}_t = \kappa {\bf 1}_{[0, \tau]} (t){\bf e}_i $, 
where 
${\bf 1}$ denotes the indicator function and 
$\{ {\bf e}_i \}_{i=1}^d$ denotes the canonical basis of ${\bf R}^d$.

For $\tau_1, \ldots, \tau_N >0$,  $i_1, \ldots, i_N \in \{1, \ldots, d\}$, and
 $\kappa_1, \ldots, \kappa_N \in \{\pm 1 \}$,
%
%we set  $\xi^{T_1, \ldots, }_{}  \in {\cal H}_{T_1+\cdots +T_m}$ by
%
 we will consider $\xi^{\tau_1, i_1, \kappa_1} * \cdots * \xi^{\tau_N, i_N, \kappa_N}
\in  {\cal H}_{\tau_1+\cdots +\tau_N}$,
whose derivative in time is given by 
$$
\sum_{k=1}^N  
\kappa_k
{\bf 1}_{[ \tau_1+\cdots +\tau_{k-1}, \tau_1+\cdots +\tau_k   ]} (t) {\bf e}_{i_k}
\qquad
\qquad
(0 \le t \le \tau_1+\cdots +\tau_N).
$$
(When $k=1$, $\tau_1+\cdots +\tau_{k-1}$ is understood to be $0$.)
For $h \in {\cal H}_{\tau}$ with $\tau >0$,  
we will write ${\cal A} h = h * \bar{h} \in   {\cal H}_{2\tau}$.

%%%%%%%%

The case $N=1$ (i.e., the elliptic case) is almost obvious,
because for any $h$ and any  $v \in {\mathbb S}^{n-1}$ there is $i$ such that
$\la v,  Q_0^{V_j} (h)  \ra  =\la v, V_i (x) \ra \neq 0$, 
which implies $\sigma_{\phi_1} (h) $ is non-degenerate for any $h$.

%%%

%%%%%%

Next, we consider the case $N=2$.
Let $\tau \in (0,T]$, where $T$ is the constant in Lemma \ref{lm.klm}.
We will prove the following:

\begin{lm}\label{lm.N2_opq}
Let $\tau$ be as above.
For any $v \in {\mathbb S}^{n-1}$, there exist $i ~(1 \le i \le d)$, $\kappa \in \{\pm 1\}$ 
such that 
$\la v,  Q_t^{V_m} (\xi^{\tau, i, \kappa}  )  \ra \neq 0$ for some $m~(1 \le m \le d)$ 
and some $t \in [0,\tau]$.
\end{lm}

\Proof
Take any $v$ and  let $W \in \Lambda$ be as in Lemma  \ref{lm.klm}.
Since we assume $N=2$, 
$W$ is of the form either $W =V_j$ or $W=[V_j, V_k]$.
If $W =V_j$ for some $j$, 
then for any $i$ and $\kappa$, 
$\la v,  Q_0^{V_j} (\xi^{\tau, i, \kappa}  )  \ra =  \la v, V_j(x) \ra \neq 0$.

Suppose that  $W=[V_j, V_k]$ for  some $j, k$.
If $\la v, V_k (x) \ra \neq 0$, then the same argument as above can still be used.
So, we may assume that $\la v, V_k (x) \ra = 0$.
Take $i=j$ and $\kappa =+1$.
By Lemma \ref{lm.abc}, {\rm (i)},  we have
\[
\Bigl|
\frac{d}{dt} \la v,  Q_t^{V_k} (\xi^{\tau, j, +1}  )  \ra 
\Bigr|
=
| \la v,  Q_t^{[V_j, V_k]} (\xi^{\tau, j, +1}  )  \ra |
\ge 
\lambda >0
\qquad
(0 \le t \le \tau).
\]
Here, we also used Lemma  \ref{lm.klm} and the choice of $W$. 
%
%the right hand side never vanishes on this interval.
%
Since the initial value of $\la v,  Q_t^{V_k} (\xi^{\tau, j, +1}  )  \ra $ is assumed to be $0$,
$\la v,  Q_t^{V_k} (\xi^{\tau, j, +1}  )  \ra \neq 0$
for any small $t>0$.
\QED

%%%%%%%%

\begin{lm}\label{lm.N2_rstu}
When $N=2$, Proposition \ref{pr.nondeg.dmc}, {\rm (i)} is true.
\end{lm}

\Proof
For $\tau \in (0,T]$,  consider 
${\cal A }\xi^{\tau, i, \kappa}= \xi^{\tau, i, \kappa} *\overline{\xi^{\tau, i, \kappa} }$
 for all $i, \kappa$
and concatenate them all, which is called $k^{\tau}$.
(The order of concatenation does not matter.)
Since there are $2d$ such $\xi^{\tau, i, \kappa}$'s, the total times length
is $2d \times 2\tau =4d\tau$.
So, $k^{\tau} \in {\cal H}_{ 4d\tau}$.
We consider ODEs (\ref{ode.def})--(\ref{ode_K.def}) driven by $k^{\tau}$.
By Lemma \ref{lm.defg}, {\rm (ii)}, 
$$
(\phi_{2\tau l} (k^{\tau}), J_{2\tau l} (k^{\tau}), J_{2\tau l} (k^{\tau})^{-1} ) 
=(x, {\rm Id}_n, {\rm Id}_n)
\qquad
\mbox{for all $l =0,1, \ldots, 2d$.}
$$
This means that, at times $2\tau, 4 \tau, \ldots, 4d\tau$,
the solution $(\phi_t, J_t, J_t^{-1})$
gets back to the initial state and starts all over again.
(If these ODEs had a drift term, this argument would fail.)

Set $\beta = 4d\tau$ and define $h^{\beta}$ for a given $h \in {\cal K}^{x, x'}$
 as follows.
On $[0, \beta]$, we set $h^{\beta}_t = k^{\tau}_t$.
On $[\beta,1]$, we set $h^{\beta}_t = h_{(t-\beta)/(1 -\beta)}$.
Then, by Lemma \ref{lm.hij},
$h^{\beta} \to h$ in ${\cal H}$ as $\tau \searrow 0$.
Moreover, by Lemma \ref{lm.defg}, {\rm (iii)}, $h^{\beta} \in {\cal K}^{x, x'}$.
By Lemma \ref{lm.N2_opq} and the way we construct $k^{\tau}$ (and $h^{\beta} $),
we have the following:
For  any $v \in {\mathbb S}^{n-1}$,  there exist $j~(1 \le j \le d)$ and
 $t \in [0,\beta]$ such that 
$\la v,  Q_t^{V_j} ( h^{\beta} )  \ra \neq 0$.
This implies non-degeneracy of $\sigma_{\phi_1} (h^{\beta} )$.
\QED

%%%%%%%%%%%%%%%

In the end of this subsection, 
we consider the case $N=3$.
If one understands the proof for this case, 
then one will easily understand the proof for the general case in the next subsection.

\begin{lm}\label{lm.N3_opq}
For sufficiently small $\tau>0$, we set $\tau_1 =\tau$ and $\tau_2 = \lambda\tau/(2M)$,
where $M>0$ is a constant given in (\ref{M.ineq}).
Then, 
for any $v \in {\mathbb S}^{n-1}$, there exist $i_1, i_2 \in \{ 1, \ldots, d\}$, 
$\kappa_1, \kappa_2 \in \{\pm 1\}$ such that 
$\la v,  Q_t^{V_m} (\xi^{\tau_1, i_1, \kappa_1} * \xi^{\tau_2, i_2, \kappa_2}  )  \ra \neq 0$
for some $m~(1 \le m \le d)$ 
and some $t \in [0, \tau_1 +\tau_2]$.
\end{lm}

\Proof
We take $\tau$  so small that $\tau_1 +\tau_2 \le T$, 
where $T$ is given in Lemma \ref{lm.klm}.
Take any $v$ and  let $W \in \Lambda$ be as in Lemma  \ref{lm.klm}.
We assume that $W$ is of the for $W = [V_j, [V_k, V_l]]$
since the other cases are easier.

On the first subinterval $[0, \tau _1]$, choose $i_1=j$.
On this interval 
$Q_t^V (\xi^{\tau_1, i_1, \kappa_1} * \xi^{\tau_2, i_2, \kappa_2} ) = Q_t^V (\xi^{\tau_1, j, \kappa_1} )$.
By Lemma \ref{lm.abc}, {\rm (i)},  we have
\begin{equation}\label{1st[].eq}
\frac{d}{dt} \la v,  Q_t^{[V_k, V_l]} (\xi^{\tau, j, \kappa_1}  )  \ra 
=
\kappa_1  \la v,  Q_t^{[V_j, [V_k, V_l]]} (\xi^{\tau, j,   \kappa_1}  )  \ra 
%\ge 
%\lambda >0
\qquad
(0 \le t \le \tau_1).
\end{equation}
On this subinterval, the right hand side of (\ref{1st[].eq}) is of constant sign,
due to Lemma \ref{lm.klm}.
If the initial value
$\la v,  Q_0^{[V_k, V_l]} (\xi^{\tau, j, \kappa_1}  ) \ra= \la v,  [V_k, V_l](x)  \ra \ge 0$,
then we choose $\kappa_1$ so that
the right hand side of (\ref{1st[].eq}) is positive.
If otherwise, 
then we choose $\kappa_1$ so that
the right hand side of (\ref{1st[].eq}) is negative.
Either way, we have $|\la v,  Q_{\tau_1}^{[V_k, V_l]} (\xi^{\tau, j, \kappa_1}  )  \ra | \ge \lambda \tau$.

On the second subinterval $[\tau _1, \tau _1+ \tau _2]$, choose $i_2=k$
and consider $\xi^{\tau_1, j, \kappa_1} * \xi^{\tau_2, k, \kappa_2}$.
By (\ref{M.ineq}) and the definition of $\tau_2$, 
\[
| \la v,  Q_{t}^{[V_k, V_l]} (\xi^{\tau, j, \kappa_1} * \xi^{\tau_2, k, \kappa_2} )  \ra 
-  
\la v,  Q_{\tau_1}^{[V_k, V_l]} (\xi^{\tau, j, \kappa_1} * \xi^{\tau_2, k, \kappa_2} )  \ra
|
\le \frac{\lambda \tau}{2}
\quad
(\tau _1 \le t \le \tau _1+ \tau _2).
\]
Hence,  we have
$|  \la v,  Q_{t}^{[V_k, V_l]} (\xi^{\tau, j, \kappa_1} * \xi^{\tau_2, k, \kappa_2} ) | \ge \lambda \tau/2$
on the second subinterval.
By Lemma \ref{lm.abc}, {\rm (i)},  we have
\begin{equation}\label{2nd[].eq}
\frac{d}{dt} \la v,  Q_t^{ V_l} (\xi^{\tau, j, \kappa_1} * \xi^{\tau_2, k, \kappa_2} )  \ra 
=
\kappa_2  \la v,  Q_t^{ [V_k, V_l]} (\xi^{\tau, j,   \kappa_1}  * \xi^{\tau_2, k, \kappa_2})  \ra 
\qquad
(\tau_1 \le t \le \tau_1+\tau_2).
\end{equation}
If the initial value 
$\la v,    Q_{\tau_1}^{ V_l} (\xi^{\tau, j, \kappa_1} * \xi^{\tau_2, k, \kappa_2} ) \ra$
of this subinterval is non-negative,
then we choose $\kappa_2$ so that
the right hand side of (\ref{2nd[].eq}) is positive.
If otherwise, 
then we choose $\kappa_2$ so that
the right hand side of (\ref{2nd[].eq}) is negative.
Either way, 
$$
|  \la v,  Q_{\tau_1+ \tau_2  }^{ V_l} (\xi^{\tau, j, \kappa_1} * \xi^{\tau_2, k, \kappa_2} )  \ra |
\ge  
\tau_2 \times \frac{\lambda \tau}{2} = \frac{(\lambda \tau)^2}{4M}. 
$$
This completes the proof.
(In fact, in order to prove this lemma 
it is enough to assume 
that $\la v, Q_{\tau_1}^{ V_l} (\xi^{\tau, j, \kappa_1} * \xi^{\tau_2, k, \kappa_2} ) \ra  = 0$ above.
However, we deliberately argued in this way for later use.)
\QED

\begin{lm}\label{lm.N3_rstu}
When $N=3$, Proposition \ref{pr.nondeg.dmc}, {\rm (i)} is true.
\end{lm}

\Proof
Let  $\tau$ be   sufficiently small and $\tau_1, \tau_2$ be as above.
Consider
\[
{\cal A} (\xi^{\tau_1, i_1, \kappa_1} * \xi^{\tau_2, i_2, \kappa_2})
\]
for all $i_1, i_2, \kappa_1, \kappa_2$.
(There are $(2d)^2$ of them.)
The concatenation of all of them is denoted by $k^{\tau}$.
(The order of concatenation does not matter.)
The total time length $\beta$ of $k^{\tau}$ is given by 
\[
\beta = (2d)^2 \cdot 2(\tau_1+\tau_2) = O(\tau)  \qquad \mbox{as $\tau \to 0$.}
\]
On $[0, \beta]$, we set $h^{\beta}_t = k^{\tau}_t$.
On $[\beta,1]$, we set $h^{\beta}_t = h_{(t-\beta)/(1 -\beta)}$.
The rest is essentially the same as the proof for the case $N=2$.
\QED

%%%%%%%%%%%%%%%%%%%%%%%%%%%%%%%%%%%%%%%%%%%%%%%%%%%%%%%%%%%%%%%%%%
%%%%%%%%%%%%%%%%%%%%%%%%%%%%%%%%%%%%%%%%%%%%%%%%%%%%%%%%%%%%%%%%%%%%%%%%
%%%%%%%%%%%%%%%%%%%%%%%%%%%%%%%%%%%%%%%%%%%%%%%%%%%%%%%%%%%%%%%%%%

\subsection{Proof of Proposition \ref{pr.nondeg.dmc} {\rm (i)}: The general case }

Now, we are in a position to
prove Proposition \ref{pr.nondeg.dmc} {\rm (i)} for the general degree $N \ge 1$.

For sufficiently small $\tau>0$, we set 
$$
\tau_1 =\tau  \quad \mbox{ and } \quad  \tau_l = 2\Bigl( \frac{\lambda\tau}{4M} \Bigr)^{2^{l-2}}
\quad
\mbox{ for $2 \le l \le N-1$.}
$$
Here, $M>0$ is a constant given in (\ref{M.ineq}).

\begin{lm}\label{lm.Ngen_opq}
Let $\tau>0$ be sufficiently small and
set $\tau_l ~(1 \le l \le N-1)$ as above.
Then, 
for any $v \in {\bf S}^{n-1}$, there exist $i_l \in \{ 1, \ldots, d\}$, 
$\kappa_l \in \{\pm 1\}$ ($1 \le l \le N-1$) such that 
$$
\la v,  Q_t^{V_m} (\xi^{\tau_1, i_1, \kappa_1} * \cdots * 
\xi^{\tau_{N-1} , i_{N-1}, \kappa_{N-1}}  )  \ra \neq 0
$$
for some $m~(1 \le m \le d)$ 
and some $t \in [0,  \tau_1 +\cdots +\tau_{N-1}  ]$.
\end{lm}

\Proof
For simplicity we write $T_l :=\tau_1 +\cdots +\tau_{l}$.
The proof is similar to the ones for Lemmas \ref{lm.N2_opq} and \ref{lm.N3_opq}.
Take any $v$ and  let $W \in \Lambda$ be as in Lemma  \ref{lm.klm}.
We assume that $W$ is of the for
 $W = [V_{j_1} ,   \cdots\cdots  [  V_{ j_{N-2} },  [V_{j_{N-1} } , V_{j_{N}} ] ]\cdots ]$
since the other cases are easier.
In this case we take
$i_l =j_l$ for $1 \le l \le N-1$
and write $\eta = \xi^{\tau_1, j_1, \kappa_1} * \cdots * 
\xi^{\tau_{N-1} , j_{N-1}, \kappa_{N-1}} $.
We will see that for a suitable choice of $\kappa_l$'s, 
$\la v,  Q_t^{V_{m} } ( \eta)  \ra \neq 0$ holds for $m=j_N$ and  $t =T_{N-1}$.

Write $U_l = [V_{j_l} ,   \cdots\cdots    [V_{j_{N-1} } , V_{j_{N}} ] \cdots ]$
for $1 \le l \le N-1$ and $U_N = V_{j_{N}}$.
On the $l$th interval ($1 \le l \le N-1$), we have 
\begin{equation}\label{3rd[].eq}
\frac{d}{dt} \la v,  Q_t^{ U_{l+1}  } (\eta)  \ra 
=
\kappa_l  \la v,  Q_t^{ U_{l}} (\eta)  \ra 
%\ge 
%\lambda >0
\qquad
\qquad
(T_{l-1} \le t \le T_l).
\end{equation}
We will prove by induction that, for suitable choices of $\kappa_i$'s, 
\begin{equation}\label{ind.kukan.ineq}
| \la v,  Q_{T_l}^{ U_{l+1}  } (\eta)  \ra|   \ge 4M \Bigl( \frac{\lambda\tau}{4M} \Bigr)^{2^{l-1}}
\quad
\mbox{ for all  $2 \le l \le N-1$.}
\end{equation}
Once this is obtained, the proof of the lemma is done since $U_N = V_{j_{N}}$.
In the same way as in Lemma \ref{lm.N3_opq},
we can prove that (\ref{ind.kukan.ineq}) holds for $l=2$ for a suitable choice of $\kappa_1$ and $\kappa_2$.
Let us assume that 
(\ref{ind.kukan.ineq}) holds up to $l-1$
for some $\kappa_1, \ldots, \kappa_{l-1}$.
By the Lipschitz continuity (\ref{M.ineq}) and the definition of $\tau_l$,
\[
|  \la v,  Q_{t }^{ U_{l}  } (\eta)  -\la v,  Q_{T_{l-1}  }^{ U_{l}  } (\eta)  \ra|  
\le 
M \tau_l
=
2M \Bigl( \frac{\lambda\tau}{4M} \Bigr)^{2^{l-2}}
\qquad
(T_{l-1} \le t \le T_l).
\]
From this estimate and (\ref{ind.kukan.ineq}) with $l-1$, 
\[
|  \la v,  Q_{t }^{ U_{l}  } (\eta) \ra | \ge 2M \Bigl( \frac{\lambda\tau}{4M} \Bigr)^{2^{l-2}}
\qquad
(T_{l-1} \le t \le T_l).
\]
Hence, the right hand side of (\ref{3rd[].eq}) is of constant sign.
If 
$\la v,  Q_{T_{l-1} }^{ U_{l+1}  } (\eta)  \ra $ is non-negative (or non-positive), 
then we 
choose 
$\kappa_l =\pm 1$ so that the right hand side of (\ref{3rd[].eq}) is positive 
  (or negative, respectively).
Then, 
it follows that 
$$
| \la v,  Q_{T_l}^{ U_{l+1}  } (\eta)  \ra|  
\ge
 2M \Bigl( \frac{\lambda\tau}{4M} \Bigr)^{2^{l-2}}  \tau_l 
=4M \Bigl( \frac{\lambda\tau}{4M} \Bigr)^{2^{l-1}},
$$
which shows that (\ref{ind.kukan.ineq}) holds up to $l$.
Thus, we have proved (\ref{ind.kukan.ineq}).
\QED

\begin{lm}\label{lm.Ngen_rstu}
Proposition \ref{pr.nondeg.dmc}, {\rm (i)} is true for any $N\ge 1$.
\end{lm}

\Proof
Once we obtain Lemma \ref{lm.Ngen_opq} above, 
the proof of the lemma is similar to that of Lemma \ref{lm.N3_rstu}. 
Let  $\tau$ be   sufficiently small and $\tau_1, \ldots, \tau_{N-1}$ be as above.
Consider
\[
{\cal A} (\xi^{\tau_1, i_1, \kappa_1} * \cdots *  \xi^{\tau_{N-1}, i_{N-1}, \kappa_{N-1}})
\]
for all $i_l, \kappa_l ~(1 \le l \le N-1)$ and concatenate them all.
(The order of concatenation does not matter.)
The total time length $\beta: = 2T_{N-1}(2d)^{N-1}$ is clearly of $O(\tau)$ as $\tau \searrow 0$. 
The rest is  the same as the proof for the case $N=3$ in Lemma \ref{lm.N3_rstu}.
\QED

%%%%%%%%%%%%%%%%%%%%%%%%%%%%%%%%%%%%%%%%%%%%%%%%%%%%%%%%%%%%%%
%%%%%%%%%%%%%%%%%%%%%%%%%%%%%%%%%%%%%%%%%%%%%%%%%%%%%%%%%%%%%
%%%%    Section    Malliavin Review
%%%%%%%%%%%%%%%%%%%%%%%%%%%%%%%%%%%%%%%%%%%%%%%%%%%%%%%%%%%%%%
%%%%%%%%%%%%%%%%%%%%%%%%%%%%%%%%%%%%%%%%%%%%%%%%%%%%%%%%%%%%%
%%%%\newpage
\section{Preliminaries}

\subsection{Preliminaries from  Malliavin calculus}

We first recall Watanabe's theory of 
generalized Wiener functionals (i.e., Watanabe distributions) in Malliavin calculus.
Most of the contents and the notations
in this section are borrowed from  Sections V.8--V.10, Ikeda and Watanabe \cite{iwbk}
with trivial modifications.
There is no new result in this section.
Shigekawa \cite{sh} and Nualart \cite{nu} are also good textbooks of Malliavin calculus.
For basic results of quasi-sure analysis, we refer to Chapter II, Malliavin \cite{ma}.

Let $({\cal W}, {\cal H}, \mu)$ be the classical Wiener space as before.  
(The results in this subsection also hold on any abstract Wiener space, however.)
The following are of particular importance in this paper:
\\
\\
{\bf (a)}~ Basics of Sobolev spaces ${\bf D}_{p,r} ({\cal K})$ of ${\cal K}$-valued 
(generalized) Wiener functionals, 
where $p \in (1, \infty)$, $r \in {\mathbb R}$, and ${\cal K}$ is a real separable Hilbert space.
As usual, we will use the spaces 
${\bf D}_{\infty} ({\cal K})= \cap_{k=1 }^{\infty} \cap_{1<p<\infty} {\bf D}_{p,k} ({\cal K})$, 
$\tilde{{\bf D}}_{\infty} ({\cal K}) 
= \cap_{k=1 }^{\infty} \cup_{1<p<\infty}  {\bf D}_{p,k} ({\cal K})$ of test functionals 
and  the spaces ${\bf D}_{-\infty} ({\cal K}) = \cup_{k=1 }^{\infty} \cup_{1<p<\infty} {\bf D}_{p,-k} ({\cal K})$, 
$\tilde{{\bf D}}_{-\infty} ({\cal K}) = \cup_{k=1 }^{\infty} \cap_{1<p<\infty} {\bf D}_{p,-k} ({\cal K})$ of 
 Watanabe distributions as in \cite{iwbk}.
When ${\cal K} ={\mathbb R}$, we simply write ${\bf D}_{p, r}$, etc.
\\
{\bf (b)}~ Meyer's equivalence of Sobolev norms. 
(Theorem 8.4, \cite{iwbk}. 
A stronger version can be found in Theorem 4.6, \cite{sh})
\\
{\bf (c)}~Pullback $T \circ F =T(F)\in \tilde{\bf D}_{-\infty}$ of tempered Schwartz distribution $T \in {\cal S}^{\prime}({\mathbb R}^n)$
on ${\mathbb R}^n$
by a non-degenerate Wiener functional $F \in {\bf D}_{\infty} ({\mathbb R}^n)$. (see Sections 5.9, \cite{iwbk}.)
\\
{\bf (d)}~A generalized version of integration by parts formula in the sense 
of Malliavin calculus
 for Watanabe distribution,
 which is given as follows (See p. 377, \cite{iwbk}):

For a non-degenerate Wiener functional
$F =(F^1, \ldots, F^n) \in {\bf D}_{\infty} ({\mathbb R}^n)$, we denote by 
$\sigma^{ij}_F (w) =  \la DF^i (w),DF^j (w)\ra_{{\cal H}}$ the $(i,j)$-component of Malliavin covariance matrix.
We denote by $\gamma^{ij}_F (w)$ the $(i,j)$-component of the inverse matrix $\sigma^{-1}_F$.
Note that $\sigma^{ij}_F \in {\bf D}_{\infty} $ and
$D \gamma^{ij}_F = \sum_{k,l} \gamma^{ik}_F ( D\sigma^{kl}_F ) \gamma^{lj}_F $.
Hence, derivatives of $\gamma^{ij}_F$ can be written in terms of
$\gamma^{ij}_F$'s and the derivatives of $\sigma^{ij}_F$'s.
Suppose 
$G \in {\bf D}_{\infty}$ and $T \in {\cal S}^{\prime} ({\mathbb R}^n)$.
Then, the following integration by parts holds;
\begin{align}
{\mathbb E} \bigl[
(\partial_i T \circ F ) \cdot G 
\bigr]
=
{\mathbb E} \bigl[
(T \circ F ) \cdot \Phi_i (\, \cdot\, ;G)
\bigr]
\label{ipb1.eq}
\end{align}
where $\Phi_i (w ;G) \in  {\bf D}_{\infty}$ is given by 
\begin{align}
\Phi_i (w ;G) &=
-\sum_{j=1}^d  
\Bigl\{
-\sum_{k,l =1}^d G(w) \gamma^{ik }_F (w)\gamma^{jl }_F (w)    \la D\sigma^{kl}_F (w),DF^j (w)\ra_{{\cal H}}
\nn\\
&
\qquad\qquad
+
\gamma^{ij }_F (w) \la DG (w),DF^j (w)\ra_{{\cal H}} + \gamma^{ij }_F (w) G (w) LF^j (w)
\Bigr\}.
\label{ipb2.eq}
\end{align}
Note that the expectations in (\ref{ipb1.eq}) are in fact  the generalized ones,
i.e.,
the pairing of $\tilde{{\bf D}}_{- \infty}$ and $\tilde{{\bf D}}_{\infty}$.
%%%%%%%%%%%%%%%%%%%%%%   %%%%%%%%%%%%%%%%%%%%%%
%%%%%%%%%%%%%%%%%%%%%%%%%%%%%%%%%%%%%%%%%%%%%
\\
\\
Let us recall Watanabe's asymptotic expansion theorem.
Let $\{ F_{\ve}\}_{0<\ve \le 1}$ be a family
of ${\mathbb R}^n$-valued Wiener functionals 
indexed by a small parameter $\ve \in (0,1]$.
If $\{ F_{\ve}\}$ admits an asymptotic expansion in ${\bf D}_{\infty} ({\mathbb R}^n)$
and their Malliavin covariance matrices are uniformly non-degenerate, 
then $T \circ  F_{\ve}$ admits an asymptotic expansion in $\tilde{\bf D}_{-\infty}$
as $\ve \searrow 0$
and each term in the expansion is obtained by formal Taylor expansion.
(Theorem 9.4, \cite{iwbk})

In this paper, however, we do not use this method.
Instead, we use a modified version of Watanabe's asymptotic expansion theorem,
which can be found in pp. 216--217, Takanobu and Watanabe \cite{tw}.

Let $\rho >0$, $\xi \in {\bf D}_{\infty}$ and $F \in {\bf D}_{\infty} ({\mathbb R}^n)$
and suppose that
\begin{equation}\label{wat1.eq}
\inf_{v \in {\mathbb S}^{n-1}}    v^* \sigma_F v  \ge \rho
\qquad
\mbox{on  \quad$\{w \in {\cal W} ~|~ |\xi (w)| \le 2 \}$.}
\end{equation}
%
%where $\sigma_F$ denotes the Malliavin covariance matrix of $F$.
%
Let $\chi:{\mathbb R} \to {\mathbb R}$ be a smooth function 
whose support is contained in $[-1, 1]$.
Then, the following proposition holds (Proposition 6.1, \cite{tw}).

\begin{pr}\label{pr.comp1}
Assume (\ref{wat1.eq}).
For every $T \in {\cal S}^{\prime}({\mathbb R}^n)$, 
$\chi(\xi) \cdot T \circ F =\chi(\xi) \cdot T (F) \in  \tilde{\bf D}_{-\infty}$
can be defined in a unique way so that the following properties hold:
\\
\noindent
{\rm (i)}~If $T_k \to T \in {\cal S}^{\prime}({\mathbb R}^n)$ as $k \to \infty$, 
then $\chi(\xi) \cdot T_k (F) \to \chi(\xi) \cdot T (F) \in \tilde{\bf D}_{-\infty}$.
\\
\noindent
{\rm (ii)}~If $T$ is given by $g \in {\cal S}({\mathbb R}^n)$,
then $\chi(\xi) \cdot T (F) = \chi(\xi)  g (F) \in {\bf D}_{\infty}$.
\end{pr}

Next, we state the asymptotic expansion theorem, which is Proposition 6.2, \cite{tw}.
Let $\{ F_{\ve}\}_{0<\ve \le 1} \subset {\bf D}_{\infty} ({\mathbb R}^n)$
and $\{ \xi_{\ve}\}_{0<\ve \le 1} \subset {\bf D}_{\infty}$
 be families
of Wiener functionals such that the following asymptotics hold:
\begin{eqnarray}
 F_{\ve}  
 &\sim&  
 f_0 +\ve f_1 + \ve^2 f_2 +\cdots 
 \qquad
 \qquad \mbox{in ${\bf D}_{\infty} ({\mathbb R}^n)$ as $\ve \searrow 0$,}
  \label{wat2.eq}
  \\
 \xi_{\ve}  
 &\sim&  
 a_0 +\ve a_1 + \ve^2 a_2 +\cdots 
 \qquad
 \qquad \mbox{in ${\bf D}_{\infty} $ as $\ve \searrow 0$.}
  \label{wat3.eq}
 \end{eqnarray}

\begin{pr}\label{pr.comp2}
Assume (\ref{wat2.eq}), (\ref{wat3.eq}) and $|a_0| \le 1/8$.
Moreover, assume that there exists $\rho >0$ independent of $\ve$ such that 
(\ref{wat1.eq}) with $F = F_{\ve}$ and $\xi = \xi_{\ve}$ holds for any $\ve \in (0,1]$.
Let $\chi:{\mathbb R} \to {\mathbb R}$ be a smooth function 
whose support is contained in $[-1, 1]$ such that $\chi (x) =1$ if $|x| \le 1/2$.
Then, we have the following asymptotic expansion:
\begin{eqnarray}
 \chi(\xi_{\ve}) \cdot T (F_{\ve}) \sim 
  \Phi_0 +\ve \Phi_1 + \ve^2 \Phi_2 +\cdots 
   \qquad
 \qquad \mbox{in $\tilde{\bf D}_{-\infty}$ as $\ve \searrow 0$.}
 \nn
 \end{eqnarray}
\end{pr}
In the above proposition, $\Phi_k \in \tilde{\bf D}_{-\infty}$ can be written as 
the $k$th coefficient of the formal Taylor expansion of $T(f_0 + [ \ve f_1 + \ve^2 f_2 +\cdots ])$.
In particular, $\Phi_0 = T(f_0)$. 
(In this paper, however, we do not need the expansion up to high order.)

%%%%%%%%%%%%%%%%%%%%%%%%%%%%%%%%%%%%%%%%%%%%%%%%%%%%%%%%%%%%%%
%%%%%%%%%%%%%%%%%%%%%%%%%%%%%%%%%%%%%%%%%%%%%%%%%%%%%%%%%%%%%%
%%%   SECTION      lower estimate
%%%%%%%%%%%%%%%%%%%%%%%%%%%%%%%%%%%%%%%%%%%%%%%%%%%%%%%%%%%%%%
%\newpage
%%%%%%%%%%%%%%%%%%%%%%%%%%%%%%%%%%%%%%%%%%%%%%%%%%%%%%%%%%%%%%
\subsection{Preliminaries from rough path theory}

In this subsection we recall the geometric rough path space 
with H\"older or Besov norm 
and quasi-sure property of rough path lift.
For basic properties of geometric rough path space, 
we refer to Lyons, Caruana, and L\'evy \cite{lcl},
and Friz and Victoir \cite{fvbk}.
For the geometric rough path space with Besov norm, 
we refer to Appendix A.2, \cite{fvbk}.
Quasi-sure property of rough path lift is summarized in Inahama \cite{in2}.
In this paper we basically assume $\alpha \in (1/3, 1/2)$ unless otherwise stated.
We always assume that 
Besov parameters $(\alpha, 4m)$ satisfy (\ref{eq.amam}), 
although some results presented in this subsection still hold 
under weaker assumptions on the parameters.

%%%%%%%%%%%%%%%%%%%%%%%%%%%%%%%%%%%%%%%%%%%%%%%%%%%%%%%%%%%%%%
%\vspace{20mm}
%%%%%%%%%%%%%%%%%%%%%%%%%%%%%%%%%%%%%%%%%%%%%%%%%%%%%%%%%%%%%%

We denote by  $G\Omega^H_{\alpha} ( {\mathbb R}^d) $  the geometric rough path space 
over ${\mathbb R}^d$ with $\alpha$-H\"older norm.
For $\beta \in (0,1]$,
let $C_0^{\beta-H}([0,1], {\mathbb R}^k)$
be the Banach space of all the ${\mathbb R}^k$-valued, $\beta$-H\"older continuous paths 
that start at $0$.
If $\alpha + \beta >1$, then the Young pairing 
\[
 G\Omega^H_{\alpha} ( {\mathbb R}^d) \times C_0^{\beta-H}([0,1], {\mathbb R}^k)
  \ni ({\bf w}, \lambda) \mapsto
  ({\bf w}, \bm{\lambda}) 
     \in  G\Omega^H_{\alpha} ( {\mathbb R}^{d+k})
     \]
is a well-defined, locally Lipschitz continuous map.
(See Section 9.4, \cite{fvbk} for instance.)

%\vspace{20mm}

Now we consider a system of RDEs driven by the Young pairing
 $({\bf w}, \bm{\lambda}) \in  G\Omega^H_{\alpha} ( {\mathbb R}^{d+1})$
of ${\bf w} \in G\Omega^H_{\alpha} ( {\mathbb R}^{d})$
and 
$\lambda \in C_0^{1-H}([0,1], {\mathbb R}^1)$.
(In most cases, we will assume $\lambda_t = \mbox{const} \times t$.)
%
%
%Let $V_{i}: {\mathbb R}^n \to {\mathbb R}^n$ be a vector field 
%which satisfies {\bf (A1)} 
%($0 \le i \le d$).
%
For vector fields $V_{i}: {\mathbb R}^n \to {\mathbb R}^n$ ($0 \le i \le d$), consider
\begin{equation}\label{rde_x.def}
dx_t = \sum_{i=1}^d  V_i ( x_t) dw_t^i + V_0 ( x_t) d\lambda_t
\qquad
\qquad
\mbox{with  \quad $x_0 =x \in {\mathbb R}^n$.}
\end{equation}
The RDEs for the Jacobian process and its inverse are given as follows;
\begin{eqnarray}
dJ_t &=& \sum_{i=1}^d  \nabla V_i ( x_t) J_t dw_t^i + \nabla V_0 ( x_t) J_t d\lambda_t
\qquad
\mbox{with $J_0 ={\rm Id} \in {\rm Mat}(n,n)$,}
\label{rde_J.def}
\\
dK_t &=& - \sum_{i=1}^d K_t  \nabla V_i ( x_t)  dw_t^i  - K_t \nabla V_0 ( x_t) d\lambda_t
\qquad
\mbox{with  $K_0 ={\rm Id} \in {\rm Mat}(n,n)$.}
\label{rde_K.def}
\end{eqnarray}
Here, 
$J, K,$ and  $\nabla V_i $ are ${\rm Mat}(n,n)$-valued.
%
%In fact, $K_t =J_t^{-1}$.
%
%

Assume that $V_i$'s are of $C_b^{4}$ for a while.
Then, a global solution of (\ref{rde_x.def})--(\ref{rde_K.def}) 
exists for any ${\bf x}$ and $\lambda$ . 
Moreover, Lyons' continuity theorem holds.
(The linear growth case is complicated and will be discussed later).
% 
%However,  Bailleul's recent result  ***  implies that the continuity theorem also holds under {\bf (A1)}.
%
%
In that case, the following map is continuous:
\[
 G\Omega^H_{\alpha} ( {\mathbb R}^d) \times C_0^{1-H}([0,1], {\mathbb R}^1)
  \ni ({\bf w}, \lambda) \mapsto
  ({\bf x}, {\bf J}, {\bf K}) \in G\Omega^H_{\alpha} ( {\mathbb R}^n \oplus {\rm Mat}(n,n)^{\oplus 2}).
    \]
(The map $({\bf w}, \lambda) \mapsto  {\bf x}$ will be denoted by 
$\Phi: G\Omega^H_{\alpha} ( {\mathbb R}^d) \times C_0^{1-H}([0,1], {\mathbb R}^1) 
\to G\Omega^H_{\alpha} ( {\mathbb R}^n$).)
Recall that in Lyons' formulation of rough path theory, the initial values of the first level paths 
must be adjusted.
When $w \in C_0^{1-H}([0,1], {\mathbb R}^d)$ and ${\bf w}$ is its natural lift, 
then the path
\[
 t \mapsto ( x + {\bf x}^1_{0,t}, {\rm Id}+{\bf J}^1_{0,t}, {\rm Id}+{\bf K}^1_{0,t} )
\]
is identical to the solution of a system (\ref{rde_x.def})--(\ref{rde_K.def}) of 
ODEs understood in the Riemann-Stieltjes sense.
Recall also that
$({\rm Id}+{\bf J}^1_{0,t})^{-1}= {\rm Id}+{\bf K}^1_{0,t}$ always holds.
%
%

%%%%%%%%%%%%%%%%%%%%%%%%%%%%%%%%%%%%%%%%%%%%%%%%%%%%%%%%%%%%%%
%%%%%%%%%%%%%%%%%%%%%%%%%%%%%%%%%%%%%%%%%%%%%%%%%%%%%%%%%%%%%%

We define a continuous function
$\Gamma:  G\Omega^H_{\alpha} ( {\mathbb R}^d) \times C_0^{1-H}([0,1], {\mathbb R}^1)
\to {\rm Mat}(n,n)$ 
as follows:
Set 
$$\Gamma ({\bf w}, \lambda)  
= ({\rm Id}+{\bf J}^1_{0,t}) C ({\bf w}, \lambda) ({\rm Id}+{\bf J}^1_{0,t})^*,
$$
where
\[
C ({\bf w}, \lambda)
:=
\int_0^1  ({\rm Id}+{\bf K}^1_{0,t}) {\bf V} ( x + {\bf x}^1_{0,t})  
{\bf V}( x + {\bf x}^1_{0,t})^*  ({\rm Id}+{\bf K}^1_{0,t})^* dt
\]
with ${\bf V} :=[V_1, \ldots, V_d] \in  {\rm Mat}(n,d)$.

If $\lambda_t =t$, then $\Gamma (\ve {\bf w}, \ve^2 \lambda) = \sigma_{X_1^{\ve}} (w)$ for $\mu$-a.a.$w$, 
where ${\bf w}$ denotes the usual Brownian rough path under $\mu$
and $X_1^{\ve}$ denotes the solution of SDE (\ref{sc_sde.def}) at $t=1$.
If $\lambda_t \equiv 0$ and ${\bf h}={\cal L}(h)$ is the natural lift of $h \in {\cal H}$, 
then
$\Gamma ({\bf h}, 0) = \sigma_{\phi_1}(h)$, the deterministic Malliavin covariance matrix 
given in (\ref{dtmmalcov.def}).

\begin{re}
In this paper we will use Lyons' continuity theorem 
only with respect to $\alpha$-H\"older topology $(1/3 <\alpha <1/2)$ and 
for $C^4_b$-coefficient vector fields. 
We do not try to extend it to the case of 
unbounded coefficient vector fields or Besov topology.
\end{re}

For $(\al, 4m)$ which satisfies (\ref{eq.amam}), 
$G\Omega^B_{\alpha, 4m} ( {\mathbb R}^d) $ denotes the geometric rough path space 
over ${\mathbb R}^d$ with $(\al, 4m)$-Besov norm.
Recall that the distance on this space is given by 
\begin{align}
d({\bf w}, \hat{\bf w}) 
&= \| {\bf w}^1- \hat{\bf w}^1 \|_{\al, 4m-B} 
+\| {\bf w}^2- \hat{\bf w}^2 \|_{2\al, 2m-B}
\nn\\
&
:=
\Bigl(
\iint_{0 \le s <t \le 1}  \frac{ | {\bf w}^1_{s,t}- \hat{\bf w}^1_{s,t}|^{4m}}
{|t-s|^{1 +4m\al }} 
dsdt
\Bigr)^{1/4m}
+
\Bigl(
\iint_{0 \le s <t \le 1}  \frac{ | {\bf w}^2_{s,t}- \hat{\bf w}^2_{s,t}|^{2m}}
{|t-s|^{1 +4m\al }} 
dsdt
\Bigr)^{1/2m}.
\nn
\end{align}
By the Besov-H\"older embedding theorem for rough path spaces,
there is a continuous embedding $G\Omega^B_{\alpha, 4m} ( {\mathbb R}^d)
 \hookrightarrow G\Omega^H_{\alpha -(1/4m)} ( {\mathbb R}^d)$.
If $\al < \al' <1/2$, there is a continuous embedding
$G\Omega^H_{\alpha'} ( {\mathbb R}^d)
\hookrightarrow G\Omega^B_{\alpha, 4m} ( {\mathbb R}^d)$.
We remark that 
we will not write the first embedding explicitly.
(For example, if we write $\Phi({\bf w}, \lambda)$ for  
$({\bf w}, \lambda) \in G\Omega^B_{\alpha, 4m} ( {\mathbb R}^d) \times C_0^{1-H}([0,1], {\mathbb R}^1)$,
then it is actually the composition of the first embedding map above and $\Phi$ with respect to 
$\{\al -1/(4m)\}$-H\"older topology.)

Note also that the Young translation by $h \in {\cal H}$ works 
well on $G\Omega^B_{\alpha, 4m} ( {\mathbb R}^d)$ under (\ref{eq.amam}).
The map $({\bf w}, h) \mapsto \tau_h ({\bf w})$
is continuous from $G\Omega^B_{\alpha, 4m} ( {\mathbb R}^d) \times {\cal H}$
to $G\Omega^B_{\alpha, 4m} ( {\mathbb R}^d)$, 
where 
$ \tau_h ({\bf w})$ is the Young translation of ${\bf w}$ by $h$.
(In the proof of the continuity, the only non-trivial components are 
the "cross integrals" in the second level paths of $\tau_h ({\bf w})$,
which are computed in Lemma 5.1, \cite{in2}.)

%%%%%%%%%%%%%%%%%%%%%%%%%%%%%%%%%%%%%%%%%%%%%%%%%%%%%%%%%%%%%%
%%%%%%%%%%%%%%%%%%%%%%%%%%%%%%%%%%%%%%%%%%%%%%%%%%%%%%%%%%%%%%
%\vspace{15mm}
%%%%%%%%%%%%%%%%%%%%%%%%%%%%%%%%%%%%%%%%%%%%%%%%%%%%%%%%%%%%%%
%%%%%%%%%%%%%%%%%%%%%%%%%%%%%%%%%%%%%%%%%%%%%%%%%%%%%%%%%%%%%%

Now we discuss quasi-sure properties of rough path lift map 
${\cal L}$ from ${\cal W}$ to $G\Omega^B_{\alpha, 4m} ( {\mathbb R}^d)$.
For $k=1,2,\ldots$ and $w \in {\cal W}$,
we denote by $w(k)$ the $k$th dyadic piecewise linear approximation of $w$ 
associated with the partition $\{ l2^{-k} ~|~ 0 \le l \le 2^k\}$ of $[0,1]$.
We set 
\[
{\cal Z}_{\al, 4m} := \bigl\{ w \in {\cal W}~|~
\mbox{ $\{ {\cal L} (w(k)) \}_{k=1}^{\infty}$ is Cauchy in $G\Omega^B_{\alpha, 4m} ( {\mathbb R}^d)$} 
\bigr\}.   
\] 
We define ${\cal L}: {\cal W} \to G\Omega^B_{\alpha, 4m} ( {\mathbb R}^d)$
by ${\cal L} (w) = \lim_{m\to \infty} {\cal L} (w(k))$ if $w \in {\cal Z}_{\al, 4m}$
and 
we do not define ${\cal L} (w)$ if $w \notin {\cal Z}_{\al, 4m}$.
(We will always use this version of ${\cal L}$.)
Note that 
${\cal H}$ and $C_0^{\beta-H}([0,1], {\mathbb R}^d)$ with $\beta \in (1/2, 1]$
are subsets of ${\cal Z}_{\al, 4m}$
and the two definition of rough path lift coincide.
Under scalar multiplication and  Cameron-Martin translation,  
 ${\cal Z}_{\al, 4m}$ is invariant.
Moreover, $c {\cal L} (w) = {\cal L} (cw)$ and $\tau_h({\cal L} (w))= {\cal L} (w+h)$ 
for any $w \in {\cal Z}_{\al, 4m}$, $c \in {\bf R}$, and $h \in {\cal H}$.

It is known that ${\cal Z}_{\al, 4m}^c$ is slim, that is $(p,r)$-capacity of 
this set is zero for any $p \in (1,\infty)$ and $r \in {\mathbb N}$.
(See Aida \cite{ai}, Inahama \cite{in1, in2}).
Therefore, from a viewpoint of quasi-sure analysis, 
the lift map ${\cal L}$ is well-defined.
(Quasi-sure property of the lift map is recently extended to the case of 
a certain class of Gaussian processes by Boediharjo, Geng, and Qian \cite{bgq}.)
Moreover, 
the map ${\cal W} \ni w \mapsto {\cal L} (w) \in G\Omega^B_{\alpha, 4m} ( {\mathbb R}^d) $
is $\infty$-quasi-continuous
(Aida \cite{ai}).
We will sometimes write ${\bf W} := {\cal L} (w)$ when it is regarded 
as a rough path space-valued random variable defined on ${\cal W}$.
Due to Lyons' continuity theorem and uniqueness of quasi-continuous modification, 
$\tilde{X}^{\ve} (\,\cdot\, , x,w) = x + \Phi(\ve {\cal L}(w), \ve^2\lambda)^1$ holds quasi-surely,
if $V_i~(0 \le i \le d)$ is of $C^3_b$.
(Here, $\lambda_t = t$.)

%%%%%%%%%%%%%%%%%%%%%%%%%%%%%%%%%%%%%%%%%%%%%%%%%%%%%%%%%%%%%%
%%%%%%%%%%%%%%%%%%%%%%%%%%%%%%%%%%%%%%%%%%%%%%%%%%%%%%%%%%%%%%
%\vspace{5mm}
%%%%%%%%%%%%%%%%%%%%%%%%%%%%%%%%%%%%%%%%%%%%%%%%%%%%%%%%%%%%%%
%%%%%%%%%%%%%%%%%%%%%%%%%%%%%%%%%%%%%%%%%%%%%%%%%%%%%%%%%%%%%

Before closing this subsection, we give a brief remark for the coefficient vector fields 
with linear growth.
Below we assume that $V_i~(0 \le i \le d)$ satisfies {\bf (A1)} and has linear growth.
In this case it is not easy to prove the existence of a global solution of RDE (\ref{rde_x.def})
for any $({\bf w}, \lambda)\in G\Omega^H_{\alpha} ( {\mathbb R}^d) \times C_0^{1-H}([0,1], {\mathbb R}^1)$.
(It must be unique if it exists.)
Hence, the Lyons-It\^o map $\Phi$ may not be defined  
on the whole space $G\Omega^H_{\alpha} ( {\mathbb R}^d) \times C_0^{1-H}([0,1], {\mathbb R}^1)$.

However,  if a global solution ${\bf x}$ exists for $({\bf w}, \lambda)$,
then we can prove with a cut-off technique that 
 a global solution exists for $({\bf w}', \lambda')$ sufficiently near $({\bf w}, \lambda)$, too.
Hence, $\Phi$ can be defined on an open subset of 
$G\Omega^H_{\alpha} ( {\mathbb R}^d) \times C_0^{1-H}([0,1], {\mathbb R}^1)$
and is continuous on it.
Let $O$ be the largest open subset with such a property.
Then, ${\cal L}({\cal H}) \times \{0\} \subset O$.

By Wong-Zakai's approximation theorem 
(a.s. convergence with respect to the sup-norm will do), 
we can see  that, for each $\ve \in (0,1]$,
$X^{\ve} (\,\cdot\, , x,w) = x + \Phi(\ve {\cal L}(w), \ve^2 \lambda)^1$,
$\mu$-a.s. on $\{ w \in {\cal W}~|~  (\ve {\cal L}(w), \ve^2 \lambda) \in O \}$.
Note that the same remark goes even if $G\Omega^H_{\alpha} ( {\mathbb R}^d)$
is replaced by $G\Omega^B_{\alpha, 4m} ( {\mathbb R}^d)$.
(In this paragraph, $\lambda_t = t$.)

%%%%%%%%%%%%%%%%%%%%%%%%%%%%%%%%%%%%%%%%%%%%%%%%%%%%%%%%%%%%%%
%%%%%%%%%%%%%%%%%%%%%%%%%%%%%%%%%%%%%%%%%%%%%%%%%%%%%%%%%%%%%%
%%%   SECTION      lower estimate
%%%%%%%%%%%%%%%%%%%%%%%%%%%%%%%%%%%%%%%%%%%%%%%%%%%%%%%%%%%%%%
%\newpage
%%%%%%%%%%%%%%%%%%%%%%%%%%%%%%%%%%%%%%%%%%%%%%%%%%%%%%%%%%%%%%
\section{Lower estimate}
The aim of this section is to prove the lower estimate
in our main theorem (Theorem \ref{tm.main}, {\rm (i)}).
The proof here is more difficult than the one for the elliptic case 
in the author's previous paper \cite{in2}.
The keys of the proof are Propositions \ref{pr.nondeg.dmc} and \ref{pr.comp2}.

In what follows,
$\Phi$ stands for the Lyons-It\^o map associated with 
the vector fields $\{V_1, \ldots, V_d; V_0\}$ with respect to 
$\{\alpha- 1/(4m)\}$-H\"older topology.
We write $\lambda^{\ve}_t =\ve^2 t$.
If ${\bf h} :={\cal L}(h)$ is the natural 
lift of $h \in {\cal H}$, then $x+\Phi ({\bf h},  0)^1_{0,t}= \phi_t(h)$,
where $\phi$ is defined by (\ref{ode.def})
and $0$ is the constant one-dimensional path $0$.

Let $X^{\ve}$ be as in (\ref{sc_sde.def}).
It is known that,
for any $h \in {\cal H}$, $X^{\ve}(1, x, w+ (h/\ve))$
admits an asymptotic expansion in ${\bf D}_{\infty}$-topology as $\ve \searrow 0$;
\[
X^{\ve}(1, x, w+ \frac{h}{\ve}) \sim f_0(h) + \ve f_1(w;h) 
+ \ve^2 f_2(w;h) +\cdots
\qquad
\mbox{in ${\bf D}_{\infty} ({\mathbb R}^n)$.}
\]
Here, $f_0(h) = \phi_1(h)=\phi(1, x, h)$
and
$ f_1(w;h) = D\phi_1 (h)\la w \ra$, which is continuous, linear in $w$.
(See Section 5, \cite{tw}).
We do not need precise information of $f_i~(i \ge 2)$ in this paper.
Obviously,  $Y^{\ve}(1, x, w+ (h/\ve))$
also admits an asymptotic expansion;
\[
Y^{\ve}(1, x, w+ \frac{h}{\ve}) \sim g_0(h) + \ve g_1(w;h) 
+ \ve^2 g_2(w;h) +\cdots
\qquad
\mbox{in ${\bf D}_{\infty} ({\mathbb R}^l)$.}
\]
Here, we set $g_i=\Pi_{{\cal V}} f_i$. In particular,  $g_0(h) = \psi_1(h)=\psi(1, x, h)$.
Hence, $g_0(h)=a$ if $h \in {\cal M}^{x, a}$.
Note that $f_1(w;h)$ and $g_1(w;h)$ are mean-zero Gaussian random vectors 
whose covariance matrices are 
$\sigma_{\phi_1}( h)$ and $\sigma_{\psi_1}( h)$, respectively.

%%%%%%%%%%%%%%%%
%\vspace{10mm}
%%%%%%%%%%%%%%%%

Let $U \subset G\Omega^B_{\al, 4m} ({\mathbb R}^d)$ be open
and ${\cal L}$ be the rough path lift map.
It suffices to show that 
\begin{equation}\label{low.est}
\liminf_{\ve \searrow 0} \ve^2 \log
\mu^{\ve}_{x, a}(U)
\ge
- \frac12 \|h\|^2_{{\cal H}}
\end{equation}
for any $h \in {\cal M}^{x, a}$ such that ${\bf h} \in U$.
Here, $\mu^{\ve}_{x, a}$ is the push-forward by $\ve {\cal L}: {\cal W} \to G\Omega^B_{\al, 4m} ({\mathbb R}^d)$
of the finite Borel measure $\theta^{\ve}_{x, a}$, 
where
$\theta^{\ve}_{x, a}$ corresponds to the positive Watanabe distribution $\delta_a ( Y^{\ve}_1)$.
%
%by $\ve {\cal L}: {\cal W} \to G\Omega^B_{\al, 4m} ({\mathbb R}^d)$.
%
Moreover, 
due to Proposition \ref{pr.nondeg.dmc},
it suffices to show (\ref{low.est}) under the additional conditions on $h$, namely, 
{\rm (i)}~$\sigma_{\phi_1}( h)$ is non-degenerate and {\rm (ii)}~ $\la h, \,\cdot\, \ra$ 
extends to a continuous linear functional on ${\cal W}$.

%%%%%%%%%%%%%%%%
%\vspace{10mm}
%%%%%%%%%%%%%%%%

For $R>0$, we set 
$$
\hat{B}_{R}=\{ {\bf w} \in G\Omega^B_{  \al, 4m} ({\mathbb R}^{d}) 
~|~ 
\|{\bf w}^1 \|_{ \al, 4m -B} +  \|{\bf w}^2 \|_{2 \al, 2m -B}^{1/2} < R
 \}
$$
and set
$\hat{B}_{R} ({\bf h})= \tau_h ( \hat{B}_{R} )$, 
where $\tau_h$ is the Young translation by $h$ on $G\Omega^B_{\al, 4m} ({\mathbb R}^{d})$.
Since $\tau_h$ is a homeomorphism for any $h \in {\cal H}$,
$\{ \hat{B}_{R} ({\bf h}) ~|~ R>0\}$ forms a fundamental system of open neighborhood around ${\bf h}$.
Since $U$ is open, 
$
\hat{B}_{R} ({\bf h}) \subset U
$
for sufficiently small $R>0$.
We will estimate the weight of $\hat{B}_{R}  ({\bf h})$ from below instead of that of $U$.

%%%%%%%%%%%%%%%%
%\vspace{10mm}
%%%%%%%%%%%%%%%%

Let $\rho >0$ be such that 
$v^* \Gamma ({\bf h}, 0) v   = v^* \sigma_{\phi_1(h)} v   \ge 2\rho$ for any $v \in {\mathbb S}^{n-1}$.
Since $\Gamma$ is continuous from 
$G\Omega^B_{\al, 4m} ({\mathbb R}^{d})\times C_0^{1-H}([0,1], {\mathbb R}^1)$,
there exist 
$R_0 >0$ and $\ve_0 \in (0,1]$ such that
\[
v^* \Gamma (\tau_h({\bf w}) , \lambda^{\ve}) v  \ge \rho
\qquad
\quad
(v \in {\mathbb S}^{n-1}, {\bf w} \in \hat{B}_{R_0}, \ve \in [0, \ve_0]).
\]
Even when $V_i$ is of linear growth, the left hand side is well-defined 
for small enough $R>0$.
Note that  the Malliavin covariance matrix of 
$$\ve^{-1} X^{\ve} (1, x, w+(h/\ve)) 
= \ve^{-1} \{ x + \Phi ( \tau_h(\ve {\bf W}) , \lambda^{\ve})^1  \}$$
equals
$\Gamma (\tau_h(\ve {\bf W}) , \lambda^{\ve})$ a.s.
Recall that 
the smallest eigenvalue of the Malliavin covariance matrix of $\ve^{-1} Y^{\ve} (1, x, w+(h/\ve))$
is larger than or equal to 
that of $\ve^{-1} X^{\ve} (1, x, w+(h/\ve))$.
Hence, 
the smallest eigenvalue of 
of the Malliavin covariance matrix of $\ve^{-1} Y^{\ve} (1, x, w+(h/\ve)) 
\ge \rho$, 
provided that $\ve{\bf W} \in \hat{B}_{R_0}$ and $\ve \in [0, \ve_0]$.

%%%%%%%%%%%%%%%%
%\vspace{10mm}
%%%%%%%%%%%%%%%%

By Cameron-Martin formula, it holds that,
for any $F \in {\bf D}_{\infty}$, 
\begin{align}
{\mathbb E} [ F  \delta_{a} (Y^{\ve}_1) ]
&=
{\mathbb E} [ \exp (- \frac{ \la h, w\ra }{\ve} - \frac{ \| h\|^2_{{\cal H}} }{2\ve^2} ) 
 F( w +\frac{h}{\ve})  \delta_{a} \Bigl(  Y^{\ve}(1, x, w+ \frac{h}{\ve}) \Bigr) ]
\nn\\
&=
e^{ - \| h\|^2_{{\cal H}}/2\ve^2 }  \ve^{-l} 
{\mathbb E} [ e^{ - \la h, w\ra /\ve } 
 F( w+\frac{h}{\ve})  
 \delta_{0} \Bigl(  \ve^{-1}[Y^{\ve}(1, x, w+ \frac{h}{\ve}) -a ]\Bigr) ].
\nn
\end{align}
Here, we used the fact that $\delta_0 (\ve\,\cdot\, ) = \ve^{-l} \delta_0 (\,\cdot\,)$.

%%%%%%%%%%%%%%%%
%\vspace{10mm}
%%%%%%%%%%%%%%%%

Let $\chi:{\mathbb R} \to {\mathbb R}$ be as in Proposition \ref{pr.comp2}.
Moreover, we assume that $\chi$ is non-increasing on $[0,\infty)$
so that $\chi$ takes values in $[0,1]$.
For sufficiently small $R>0$,
we have 
\begin{align}
\mu^{\ve}_{x, a}( \hat{B}_R ({\bf h}) ) 
&=
\int  I_{\hat{B}_R ({\bf h})} ({\bf w})  \mu^{\ve}_{x, a} (d {\bf w})
=
\int  I_{\hat{B}_R } (\tau_{-h} ({\bf w}))  \mu^{\ve}_{x, a} (d {\bf w})
\nn
\\
&=
\int  I_{\hat{B}_R } (\tau_{-h} (\ve {\bf W}))  \theta^{\ve}_{x, a} (dw)
\nn
\\
&\ge
\int  \chi \Bigl( 
\frac{ \|\tau_{-h} (\ve {\bf W})^1 \|_{\alpha, 4m -B}^{4m} 
+ 
\|\tau_{-h} (\ve {\bf W})^2 \|_{2\alpha, 2m -B}^{2m}}{R^{4m}}
 \Bigr)  \theta^{\ve}_{x, a} (dw)
\nn
\\
&= 
{\mathbb E}  \Bigl[   
\chi \Bigl( 
\frac{ \|\tau_{-h} (\ve {\bf W})^1 \|_{\alpha, 4m -B}^{4m} 
+ 
\|\tau_{-h} (\ve {\bf W})^2 \|_{2\alpha, 2m -B}^{2m}}{R^{4m}}
 \Bigr) 
 \delta_{a} (Y^{\ve}_1) 
 \Bigr]
\nn\\
&=
e^{ - \| h\|^2_{{\cal H}}/2\ve^2 }  \ve^{-l} 
{\mathbb E} \Bigl[ e^{ - \la h, w\ra /\ve } 
  \chi \Bigl( 
\frac{ \| (\ve {\bf W})^1 \|_{\alpha, 4m -B}^{4m} 
+ 
\| (\ve {\bf W})^2 \|_{2\alpha, 2m -B}^{2m}}{R^{4m}}
 \Bigr) 
 \nn
 \\
 & \qquad\qquad \quad \qquad\qquad \quad
  \times \delta_{0} \Bigl(  \ve^{-1}[Y^{\ve}(1, x, w+ \frac{h}{\ve}) -a ]\Bigr) \Bigr].
\nn
%
%{\mathbb E} [ F  \delta_{a} (Y^{\ve}_1) ]
%&=
\end{align}
Note that $w ={\bf W}^1$ and 
if $\| (\ve {\bf W})^1 \|_{\alpha, 4m -B} \le R$, then $ \|w\|_{\infty} \le c_1 R/\ve$,
where 
$c_1>0$ is the operator norm of the embedding of 
the usual path space with $(\alpha, 4m)$-Besov norm into the one with the sup-norm.
Therefore, 
$e^{ - \la h, w\ra /\ve }  \ge \exp ( - c_2 R/\ve^2)$,
where $c_2 := c_1 \| \la h, \,\cdot\,\ra \|_{{\cal W}^*}$ is a positive constant 
independent of $R, \ve$.
Noting that the positive Watanabe distribution 
 $\delta_{0} (  \ve^{-1}[Y^{\ve}(1, x, w+ \frac{h}{\ve}) -a ])$
is in fact a finite measure on ${\cal W}$ by Sugita's theorem \cite{su}, 
we see that
\begin{align}
\mu^{\ve}_{x, a}( \hat{B}_R ({\bf h}) ) 
&\ge
e^{ - \| h\|^2_{{\cal H}}/2\ve^2 }  \ve^{-l}  e^{ - c_2 R /\ve^2 } \times
\nn
\\
&\quad
{\mathbb E} \Bigl[
  \chi \Bigl( 
\frac{ \| (\ve {\bf W})^1 \|_{\alpha, 4m -B}^{4m} 
+ 
\| (\ve {\bf W})^2 \|_{2\alpha, 2m -B}^{2m}}{R^{4m}}
 \Bigr) 
% \nn
% \\
% & \qquad\qquad \quad \qquad\qquad \quad
 \delta_{0} \Bigl(  \ve^{-1}[Y^{\ve}(1, x, w+ \frac{h}{\ve}) -a ]\Bigr) \Bigr].
\nn
\end{align}

Assume that the logarithm
of the generalized expectation on the right hand side above 
converges to some (finite) real quantity  (for each $R>0$ small enough)
as $\ve \searrow 0$, which we will prove later.
Then, we have 
\[
\liminf_{\ve \searrow 0}
\ve^2 
\log
\mu^{\ve}_{x, a}( \hat{B}_{R_0} ({\bf h}) )  
\ge
\liminf_{\ve \searrow 0}
\ve^2 
\log
\mu^{\ve}_{x, a}( \hat{B}_{R} ({\bf h}) )  
\ge
- \frac12 \| h\|^2_{{\cal H}} - c_2 R
\]
if $R \in (0, R_0)$.
Letting $R \searrow 0$, we obtain the desired estimate (\ref{low.est}).

%%%%%%%%%%%%%%%%
%\vspace{10mm}
%%%%%%%%%%%%%%%%

Now we use Proposition \ref{pr.comp2} with 
$\xi_{\ve} = \{ \| (\ve {\bf W})^1 \|_{\alpha, 4m -B}^{4m} 
+ 
\| (\ve {\bf W})^2 \|_{2\alpha, 2m -B}^{2m} \}/R^{4m}$,
$T =\delta_0$,
and 
$F_{\ve} =  \ve^{-1}[Y^{\ve}(1, x, w+ (h/\ve)) -a ]$.
Notice that if $0< R \le 2^{-1/m} R_0$, then 
the condition (\ref{wat1.eq}) is satisfied with $F=F_{\ve}$, $\xi =\xi_{\ve}$, 
and $\rho >0$ defined as above.
Proposition \ref{pr.comp2} implies that
\[
\chi(\xi_{\ve}) \cdot \delta_0 (F_{\ve})
=
\delta_0 (g_1(w;h))
+O(\ve)
\qquad
\mbox{in $\tilde{{\bf D}}_{- \infty}$ as $\ve \searrow 0$}
\]
Since $g_1(w;h)$ is a non-degenerate Gaussian random variable 
taking values in ${\mathbb R}^l$,
its law has a strictly positive density.
Hence, we have
$$
\lim_{\ve \searrow 0} 
{\mathbb E} [\chi(\xi_{\ve}) \cdot \delta_0 (F_{\ve})]
=
{\mathbb E} [\delta_0 (g_1(w;h)) ] \in (0, \infty).
$$
Therefore, 
$\lim_{\ve \searrow 0} \log {\mathbb E} [\chi(\xi_{\ve}) \cdot \delta_0 (F_{\ve})] \in (- \infty, \infty)$,
which completes the proof
of the lower estimate of our main theorem.

\begin{re}\label{re.added}
In the last part of the proof above, we used a modifed version of 
Watanabe's asymptotic expansion (Proposition \ref{pr.comp2}).
However, as in the proof of the elliptic case in \cite{in2},
it may also be nice to use the standard version (Theorem 9.4, \cite{iwbk})
after proving the uniform non-degeneracy of $\ve^{-1}[Y^{\ve}(1, x, w+ (h/\ve)) -a ]$
when the deterministic Malliavin covariance matrix at $h$ is non-degenerate. 
(Loosely speaking, the authors of \cite{bmn} argue in this way, for instance.)   
   
In the hypoeliptic case, however, the proof of 
  uniform non-degeneracy becomes more difficult.
We need to combine {\rm (i)} Kusuoka-Stroock's bound (\ref{ks_malcov.ineq}) below for $Y^{\ve}$ instead of  $X^{\ve}$
and  {\rm (ii)} the Schilder-type LDP for Brownian rough path.
  \end{re}

%%%%%%%%%%%%%%%%%%%%%%%%%%%%%%%%%%%%%%%%%%%%%%%%%%%%%%%%%%%%%%
%%%%%%%%%%%%%%%%%%%%%%%%%%%%%%%%%%%%%%%%%%%%%%%%%%%%%%%%%%%%%%
%%%   SECTION      lower estimate
%%%%%%%%%%%%%%%%%%%%%%%%%%%%%%%%%%%%%%%%%%%%%%%%%%%%%%%%%%%%%%
%\newpage
%%%%%%%%%%%%%%%%%%%%%%%%%%%%%%%%%%%%%%%%%%%%%%%%%%%%%%%%%%%%%%
\section{Upper estimate}
The aim of this section is to prove the upper estimate
in our main theorem (Theorem \ref{tm.main}, {\rm (i)}).
The proof here is  similar to the one for the elliptic case 
in \cite{in2} and is a modification of it.

In this section, we will often use the following fact;
For $f, g: (0,1] \to [0, \infty)$, it holds that 
$
\limsup_{\ve \searrow 0} \ve^2 \log ( A_{\ve} + B_{\ve}  )  
\le
[ \limsup_{\ve \searrow 0} \ve^2 \log  A_{\ve}  ] \vee 
 [ \limsup_{\ve \searrow 0} \ve^2 \log  B_{\ve} ].
$ 
%
%
%
%We will denote by $\Phi_0$ the Lyons-It\^o map associated with 
%the vector fields $\{V_1, \ldots, V_d; V_0\}$.
%For a geometric rough path ${\bf w}$ over ${\mathbb R}^d$, 
%$({\bf w}, {\bf 0}$) stands for the Young pairing of ${\bf w}$ and the 
%constant one-dimensional path $0$.
%It is  a geometric rough path  over ${\mathbb R}^{d+1}$.
%
%
%If ${\bf h} ={\cal L}(h)$ is the natural 
%lift of $h \in {\cal H}$, then $x+\Phi ({\bf h},  0)^1_{0,t}= \phi_t(h)$,
%where $\phi$ is defined by (\ref{ode.def}), 
%where $0$ is the constant one-dimensional path $0$.
%
%
%
We assume without loss of generality that
${\cal V} = {\mathbb R}^l \times \{ {\bf 0}_{n-l}\}$ 
so that $Y^{\ve}_t = \Pi_{{\cal V}} (X^{\ve}_t) =( X^{\ve,1}_t, \ldots , X^{\ve, l}_t)$.
(This assumption is just for notational simplicity.)

%%%%%%%%%%%%%%%%%%%%
%\vspace{10mm}
%%%%%%%%%%%%%%%%%%%%%

{\bf [Step 1]}~
We divide the proof into three steps. 
The first step is to show that 
\begin{equation}\label{limsmball.ineq}
\lim_{R \searrow 0}  \limsup_{\ve \searrow 0}  
\ve^2 \log \mu_{x,a}^{\ve} ( B_{R} ({\bf w}) )  \le   - I({\bf w}),
\qquad 
{\bf w} \in G\Omega^B_{\al, 4m}  ({\mathbb R}^d),
\end{equation}
where 
\[
B_{R} ({\bf w})
=\{ {\bf v}\in G\Omega^B_{  \al, 4m} ({\mathbb R}^{d}) ~|~ \|{\bf v}^i -{\bf w}^i \|_{i \al, 4m/i -B}  <R^{i}
\quad (i=1,2) \}.
\]

%%%%%%%%%%%%%%%%%%%%%%%%%%%%%%%%%%%%%%%%%%%%%%%%%%%%%%%%%%%%%%
%%%%%%%%%%%%%%%%%%%%%%%%%%%%%%%%%%%%%%%%%%%%%%%%%%%%%%%%%%%%%%
%\vspace{10mm}
%%%   
%%%%%%%%%%%%%%%%%%%%%%%%%%%%%%%%%

First, we consider the case ${\bf w} = {\bf h}$, where $h \in {\cal H} \setminus {\cal M}^{x,a}$.
We write
$\tilde{a} := \Pi_{{\cal V}} (x+ \Phi ( {\bf h},  0)^1_{0,1} ) ~(\neq a)$.
%
%For simplicity we denote the left hand side by $\tilde{a}$.
%
%
By the arguments in the previous section, 
even if $V_i$'s admits linear growth,
there exist $\ve_0 >0$ and $R >0$ such that 
$({\bf v}, \ve) \mapsto 
\Phi ( {\bf v}, \lambda^{\ve})$ is well-defined and continuous on $B_{2^{1/4m}R} ({\bf h}) \times [0, \ve_0)$.
Moreover, 
we may assume that
$|\Pi_{{\cal V}} (x+ \Phi ( {\bf v}, \lambda^{\ve})^1_{0,1} ) - \tilde{a} | \le |a - \tilde{a}|/3$
holds for all $0 \le \ve \le \ve_0$ and ${\bf v} \in B_{2^{1/4m}R} ({\bf h})$.
Note that for a fixed $\ve \in [0, \ve_0)$,
$X^{\ve} (\,\cdot\, , x , w) = x+ \Phi ( \ve {\bf W}, \lambda^{\ve})^1$
holds $\mu$-a.s. on $\{ w \in {\cal W}~|~ \ve {\bf W} \in B_{2^{1/4m}R} ({\bf h})\}$.
%
%
%If we show $\mu_{x,a}^{\ve} ( B_{R} ({\bf h}) ) =0$ for such $\ve$ and $R$,
%then (\ref{limsmball.ineq}) immediately follows for this case.
%
%
\\
\\

Let us verify that $\mu_{x,a}^{\ve} ( B_{R} ({\bf h}) ) =0$ as follows. 
Then, (\ref{limsmball.ineq}) immediately follows.
Let $\chi: {\mathbb R} \to [0,1]$ be a smooth even function such that
$\chi = 1$ on $[0,1]$ and $\chi  = 0$ on $[2, \infty)$ and  non-increasing on $[1,2]$.
%
%
%Set $\chi (| ( Y^{\ve}_1 -a) /  \eta|^2 )$, where $\eta : =  |a - \tilde{a}|/3$.
%Then, $\delta_{a} ( Y^{\ve}_1 ) =   \chi (| ( Y^{\ve}_1 -a) /  \eta|^2 ) \cdot \delta_{a} ( Y^{\ve}_1 )$
%in $\tilde{\bf D}_{- \infty}$.
%
%
%By Sugita's theorem,  $\theta^{\ve}_{x,a} (dw)
% =  \chi (| (\Pi_{{\cal V}} (x+ \Phi ( \ve{\bf W}, \lambda^{\ve})^1_{0,1} )    -a) /  \eta|^2 ) 
%  \cdot \theta^{\ve}_{x,a} (dw) $,
%since $ \Pi_{{\cal V}} (x+ \Phi ( \ve{\bf W}, \lambda^{\ve})^1_{0,1} )$ 
%is the $\infty$-quasi-redefinition 
%of $Y^{\ve}_1 = Y^{\ve} (1, x, w)$.
%
%
%
Let $\{f_k\}$ be a sequence of continuous functions on ${\mathbb R}^l$ such that 
$f_k \to \delta_a$ in ${\cal S}^{\prime} ({\mathbb R}^l)$ as $k \to \infty$.
We may assume that the support of $f_k$ 
is contained in $\{ \xi \in {\mathbb R}^l ~|~ |\xi -a|< |a -\tilde{a}|/3 \}$ for any $k$.

Then we have, 
\begin{align}
\mu_{x,a}^{\ve} ( B_{R} ({\bf h}) ) 
&= \theta_{x,a}^{\ve} (\{ w \in {\cal W} ~|~\ve {\bf W} \in B_{R} ({\bf h}) \} )
\nn\\
&\le
\int_{{\cal W}}
\prod_{i=1}^2
\chi \bigl( \| \ve^i {\bf W}^i - {\bf h}^i \|^{4m/i}_{i\al, 4m/i -B} / R^{4m} \bigr)
\theta^{\ve}_{x,a} (dw)
\nn\\
&=
{\mathbb E}
\Bigl[
\prod_{i=1}^2
\chi \bigl( \| \ve^i {\bf W}^i - {\bf h}^i \|^{4m/i}_{i\al, 4m/i -B} / R^{4m} \bigr)
\delta_a (Y^{\ve}_1)
\Bigr]
\nn\\
&=
\lim_{k \to \infty}
{\mathbb E}\Bigl[
\prod_{i=1}^2
\chi \bigl( \| \ve^i {\bf W}^i - {\bf h}^i \|^{4m/i}_{i\al, 4m/i -B} / R^{4m} \bigr)
f_k (Y^{\ve}_1)
\Bigr]
\nn\\
&=
\lim_{k \to \infty}
{\mathbb E}\Bigl[
\prod_{i=1}^2
\chi \bigl( \| \ve^i {\bf W}^i - {\bf h}^i \|^{4m/i}_{i\al, 4m/i -B} / R^{4m} \bigr)
f_k ( \Pi_{{\cal V}}(x+ \Phi ( \ve {\bf W}, \lambda^{\ve})^1_{0,1} ))
\Bigr]
=0.
\nn
\end{align}

%%%%%%%%%%%%%%%%%%%%%%%%%%%%%%%%%%%%%%%%%%%%%%%%%%%%%%%%%%%%%%
%%%%%%%%%%%%%%%%%%%%%%%%%%%%%%%%%%%%%%%%%%%%%%%%%%%%%%%%%%%%%%
%\vspace{10mm}
%%%   
%%%%%%%%%%%%%%%%%%%%%%%%%%%%%%%%%%%%%%%%%%%%%%%%%%%%%%%%%%%%%%

Next, let us consider the other case, namely 
${\bf w} \in G\Omega^B_{\al, 4m}  ({\mathbb R}^d) \setminus {\cal L} ({\cal H} \setminus {\cal M}^{x,a})$.
Note that 
$\| D^r  X^{\ve}_1 \|_{L^p}$ is bounded in $\ve$ for any $p \in (1, \infty)$ and 
$r=0,1,2, \ldots$,
where $D$ stands for the ${\cal H}$-derivative.
%
%Note also that 
%
%
Recall  that the Malliavin covariance matrix of $X^{\ve}_1$ satisfies 
the following estimate;
\begin{equation}\label{ks_malcov.ineq}
\| (\det \sigma_{X^{\ve}_1})^{-1} \|_{L^p} 
\le
K_1(p) \ve^{-K_2} 
\qquad
\quad
\mbox{for all $p \in (1, \infty)$ and $\ve \in (0,1]$}.
\end{equation}
Here, the constant $K_1(p) >0$ may depend on $p$ but not on $\ve$ 
and 
the constant $K_2 >0$ does not depend on $p, \ve$ (the starting point $x$ is fixed here).
This can be found in Section V-10, \cite{iwbk} or originally in Kusuoka-Stroock \cite{ksII}.

%%%%%%%%%%%%%%%%%%%%%%%%%%%%%%%%%%%%%%%%%%%%%%%%%%%%%%%%%%%%%%
%%%%%%%%%%%%%%%%%%%%%%%%%%%%%%%%%%%%%%%%%%%%%%%%%%%%%%%%%%%%%%
%\vspace{10mm}
%%%   
%%%%%%%%%%%%%%%%%%%%%%%%%%%%%%%%%%%%%%%%%%%%%%%%%%%%%%%%%%%%%%

%
Note that $\| \ve {\bf W}^1 - {\bf w}^1 \|^{4m}_{\al, 4m -B}$ and 
$\| \ve^2 {\bf W}^2 - {\bf w}^2 \|^{2m}_{2\al, 2m -B}$ 
belong to the $4m$-th order inhomogeneous Wiener chaos.
Since their $L^2$-norms are bounded in $\ve$, 
so are their $(p,r)$-Sobolev norms for any $(p,r)$.

Set $G(u_1, \ldots, u_l) = \prod_{j=1}^l (u_j - a_j)^+$,
which is a continuous function from ${\mathbb R}^l$ to ${\mathbb R}$
with polynomial growth
and satisfies $\partial_1^2 \cdots \partial_l^2 G (u)=\delta_{a} (u)$ in 
the sense of Schwartz distributions on ${\mathbb R}^l$.
It is straight forward to check that 
$( \partial_1^2 \cdots \partial_l^2 G)  \circ \Pi_{{\cal V}} 
= ( \partial_1^2 \cdots \partial_l^2 )(G \circ \Pi_{{\cal V}} ) $
as a 
Schwartz distribution on ${\mathbb R}^n$.

Then, we have 
\begin{align}
\mu_{x,a}^{\ve} ( B_{R} ({\bf w}) ) 
&= 
\theta_{x,a}^{\ve} (\{ w \in {\cal W} ~|~\ve {\bf W} \in B_{R} ({\bf w}) \} )
\nn\\
&\le
{\mathbb E}
\Bigl[
\prod_{i=1}^2
\chi \bigl( \| \ve^i {\bf W}^i - {\bf w}^i \|^{4m/i}_{i\al, 4m/i -B} / R^{4m} \bigr)
\cdot
( \partial_1^2 \cdots \partial_l^2 G) (Y^{\ve}_1)
\Bigr]
\nn\\
&\le
{\mathbb E}
\Bigl[
\prod_{i=1}^2
\chi \bigl( \| \ve^i {\bf W}^i - {\bf w}^i \|^{4m/i}_{i\al, 4m/i -B} / R^{4m} \bigr)
\cdot
( \partial_1^2 \cdots \partial_l^2 )(G \circ \Pi_{{\cal V}} ) (X^{\ve}_1)
\Bigr].
\label{ups1.eq}
\end{align}

%%%%%%%%%%%%%%%%%%%%%%%%%%%%%%%%%%%%%%%%%%%%%%%%%%%%%%%%%%%%%%
%%%%%%%%%%%%%%%%%%%%%%%%%%%%%%%%%%%%%%%%%%%%%%%%%%%%%%%%%%%%%%
%\vspace{20mm}
%%%   
%%%%%%%%%%%%%%%%%%%%%%%%%%%%%%%%%%%%%%%%%%%%%%%%%%%%%%%%%%%%%%

Now we use the integration by parts formula for Watanabe distributions
as in p. 377, \cite{iwbk}.
Then, the right hand side of (\ref{ups1.eq}) is equal to 
a finite sum of the following form;
\begin{align}
\sum_{j,j'}
{\mathbb E}
\Bigl[
F_{j,j'}^{\ve} \cdot
\chi^{(j)} \Bigl( \frac{\| \ve {\bf W}^1 - {\bf w}^1 \|^{4m}_{\al, 4m -B} }{ R^{4m} }\Bigr)
&
\chi^{(j')} \Bigl( \frac{ \| \ve^2 {\bf W}^2 - {\bf w}^2 \|^{2m}_{2\al, 2m -B} }{  R^{4m} }\Bigr)
(G \circ \Pi_{{\cal V}} ) (X^{\ve}_1)
\Bigr].
\label{ups2.eq}
\end{align}
Here, $F_{j,j'}^{\ve} (w) =F_{j,j'}(\ve, w)$ is a polynomial 
in components of 
(i)~$X^{\ve}_1$ and its derivatives,
(ii)~$\| \ve^i {\bf W}^i - {\bf w}^i \|^{4m/i}_{i\al, 4m/i -B}$ and its derivatives,
(iii)~ $\sigma_{ X^{\ve}_1}$, which is a Malliavin covariance matrix of $X^{\ve}_1$,
and (iv)~ $\gamma_{ X^{\ve}_1}= (\sigma_{ X^{\ve}_1})^{-1}$.
Note that derivatives of $\gamma_{ X^{\ve}_1}$ do not appear.
It is important that the right hand side of (\ref{ups2.eq}) 
is not a generalized expectation anymore.
By (\ref{ks_malcov.ineq})
there exists a constant $K >0$ such that $F_{j,j'}^{\ve}$ is $O(\ve^{-K})$ in any $L^p$.
(Below $K$ may change from line to line. The exact value of $K$ is of no importance.)

%%%%%%%%%%%%%%%%
%\vspace{10mm}
%%%%%%%%%%%%%%%%

Take any $p, q \in (1, \infty)$ such that $1/p + 1/q =1$.
By H\"older's inequality, 
the right hand side of (\ref{ups2.eq}) is dominated by
\begin{align}
C\ve^{-K}
\mu
\Bigl( \| \ve^i {\bf W}^i - {\bf w}^i \|_{i \al, 4m/i -B}^{1/i}   \le 2^{\frac{1}{4m}  } R
~(i=1,2)
 \Bigr)^{\frac{1}{q} }
=
C\ve^{-K} \mu (\ve {\bf W} \in B_{2^{ 1/4m  }   R} ({\bf w}))^{\frac{1}{q} }.
\nn
\end{align}
Here, $C=C_{p,q} >0$ is a constant independent of $\ve$.
By the large deviation principle of Schilder-type for the scaled Brownian 
rough path $\ve {\bf W}$ on 
$G\Omega^B_{\al, 4m} ({\mathbb R}^d)$, we have
\[
\limsup_{\ve \searrow 0} \ve^2 \log \mu_{x,a}^{\ve} ( B_{R} ({\bf w}) )  
\le
- \frac{1}{q}  \inf\{ \|h\|^2_{{\cal H} } /2 ~|~ h \in {\cal H}, {\cal L}(h) \in B_{2^{1/4m} R} ({\bf w})\}.
\]
Letting $q \searrow 1$, we have
\begin{eqnarray*}
\limsup_{\ve \searrow 0} \ve^2 \log \mu_{x,a}^{\ve} ( B_{R} ({\bf w}) )  
&\le&
- \inf\{ \|h\|^2_{{\cal H} } /2 ~|~ h \in {\cal H},   {\cal L}(h) \in B_{2^{1/4m} R} ({\bf w})\}
\\
&=&
- \inf\{ I_{Sch}({\bf v}) ~|~ {\bf v} \in B_{2^{1/4m} R} ({\bf w}) \}.
\end{eqnarray*}
Since the good rate function $I_{Sch}: G\Omega^B_{\al, 4m} ({\mathbb R}^d) \to [0, \infty]$
is lower semicontinuous, 
the limit of the right hand side as $R \searrow 0$ is dominated by $- I_{Sch}({\bf w})$.
This proves (\ref{limsmball.ineq}).
(Here, $I_{Sch}( {\cal L}  (h)) := \|h\|^2_{{\cal H} } /2$ and 
$I_{Sch}({\bf w}) :=\infty$ if ${\bf w}$ is not the natural lift of any $h \in {\cal H}$.)

%%%%%%%%%%%%%%%%
%\vspace{10mm}

%
%%%%%%%%%%%%%%%%%%%%%
{\bf [Step 2]}~
The second step is to prove the upper bound in our main theorem
 (Theorem \ref{tm.main}, {\rm (i)})
when $A$ is a compact set in $G\Omega^B_{\al, 4m} ({\mathbb R}^d)$.
Let $N \in {\mathbb N}$ be sufficiently large.
For any ${\bf w} \in A$,  take $R =R_{N, {\bf w}} >0$ small enough so that 
\[
\limsup_{\ve \searrow 0}  \ve^2 \log \mu_{x,a}^{\ve} ( B_R ({\bf w}) ) \le 
\begin{cases}
-N & ( \mbox{if  $I({\bf w})= \infty$}), \\
- I({\bf w}) +N^{-1}  &   ( \mbox{if  $I({\bf w}) < \infty$}).
\end{cases}
\]
%
% 
%$\limsup_{\ve \searrow 0}  \ve^2 \log \mu_{a,a'}^{\ve} ( B_r (Y) ) \le - I(Y) +N^{-1}$ if $I(Y) < \infty$.
%
%
The union of such open balls over ${\bf w} \in A$ covers the compact set $A$.
Hence, there are finitely many ${\bf w}_1, \ldots, {\bf w}_k \in A$ such that 
$A \subset  \cup_{i=1}^k B_{R_i} ({\bf w}_i)$, where $R_i = R(N, {\bf w}_i)$.
By using the remark in the beginning of this section,  we see that
\begin{align}
\limsup_{\ve \searrow 0}  \ve^2 \log \mu_{x,a}^{\ve} ( A ) 
&\le
(-N ) \vee \max\{ -I( {\bf w}_i) +N^{-1}  ~|~ 1 \le i  \le k, \,\,  I({\bf w}_i) < \infty  \}
\nn\\
&\le 
(-N) \vee \bigl(  - \inf_{h \in {\cal L}^{-1} (A) \cap {\cal M}^{x,a}}  \| h\|^2_{{\cal H}}/2   +N^{-1} \bigr).
\nn
\end{align}
Letting $N \to \infty$, we obtain 
\[
\limsup_{\ve \searrow 0}  \ve^2 \log \mu_{x,a}^{\ve} ( A ) 
\le
- \inf
\{
\| h\|^2_{{\cal H}}/2
~|~ 
h \in {\cal M}^{x,a}, \, {\cal L} (h) \in A
\}.
\]
Thus, we have obtained the upper estimate for the compact case.

%%%%%%%%%%%%%%%%
%\vspace{10mm}

%
%%%%%%%%%%%%%%%%%%%%%
{\bf [Step 3]}~
In this final step we will prove the upper bound in our main theorem
 (Theorem \ref{tm.main}, {\rm (i)})
when $A$ is a closed set in $G\Omega^B_{\al, 4m} ({\mathbb R}^d)$.
%
%

%For sufficiently large $\rho >0$, consider 
%$A = \{ A \cap \overline{B_{\rho}({\bf 0}) } \} \cup \{ A \cap (\overline{B_{\rho}({\bf 0}) })^c\}$,
%where ${\bf 0}$ is the trivial rough path.

Take $\alpha'$ slightly greater than $\alpha$ so that the condition (\ref{eq.amam})
still holds for $(\alpha', m)$.
Then, the continuous embedding 
$G\Omega^B_{\al', 4m} ({\mathbb R}^d) \hookrightarrow G\Omega^B_{\al, 4m} ({\mathbb R}^d)$
is in fact compact, 
which means that any bounded set in $G\Omega^B_{\al', 4m} ({\mathbb R}^d)$
is precompact in $G\Omega^B_{\al, 4m} ({\mathbb R}^d)$.
(See \cite{in2}.)

For $\rho >0$, denote by $B^{\prime}_{\rho}({\bf 0})$ the ball with respect to 
$(\al', 4m)$-Besov norm of radius $\rho$
and centered at the trivial rough path ${\bf 0}$.
Then, $B^{\prime}_{\rho}({\bf 0})$ is precompact 
with respect to 
$(\al, 4m)$-Besov topology.

%$A = \{ A \cap \overline{B_{\rho}({\bf 0}) } \} \cup \{ A \cap (\overline{B_{\rho}({\bf 0}) })^c\}$,
%where ${\bf 0}$ is the trivial rough path.

%%%%%%%%%%%%%%%%
%\vspace{10mm}
%%%%%%%%%%%%%%%%

Then, $\overline{A \cap B^{\prime}_{\rho}({\bf 0})}$ is  compact with respect to $(\alpha, 4m)$-Besov topology
and is included in $A =\bar{A}$, 
where the closure is taken with respect to 
$(\al, 4m)$-Besov topology.
Hence, we can use the argument in the previous step:
\begin{align}
\limsup_{\ve \searrow 0}  \ve^2 \log \mu_{x,a}^{\ve} ( \overline{A \cap B^{\prime}_{\rho}({\bf 0})} ) 
&\le
- \inf
\{
\frac{  \| h\|^2_{{\cal H}}  }{ 2}
~|~ 
h \in {\cal M}^{x,a}, \, {\cal L} (h) \in  \overline{A \cap B^{\prime}_{\rho}({\bf 0})}  
\}
\nn\\
&\le 
- \inf
\{
\frac{  \| h\|^2_{{\cal H}}  }{ 2}
~|~ 
h \in {\cal M}^{x,a}, \, {\cal L} (h) \in  A   
\}.
\nn
\end{align}

On the other hand, the weight of $B^{\prime}_{\rho}({\bf 0}) ^c$
is dominated as follows.
\begin{eqnarray}
\mu_{x,a}^{\ve} ( B^{\prime}_{\rho}({\bf 0})^c)
&=&
\theta_{x,a}^{\ve} (\{w \in {\cal W} ~|~ \ve {\bf W}  \notin B^{\prime}_{\rho}({\bf 0})  \})
\nn\\
&\le&
\|  \delta_a (Y^{\ve}_1)\|_{2,-r} 
{\rm Cap}_{2,r}(\{w \in {\cal W} ~|~ \ve {\bf W}  \notin B^{\prime}_{\rho}({\bf 0})  \})
\nn\\
&\le&
\|  \delta_a (Y^{\ve}_1)\|_{2,-r} \Bigl[
{\rm Cap}_{2,r}(\{w \in {\cal W} 
~|~
\| {\bf W}^1 \|_{\al', 4m -B} \ge \frac{\rho}{\ve}  \})
\nn\\
&&  \qquad\qquad\qquad
+
{\rm Cap}_{2,r}(\{w \in {\cal W} 
~|~
\| {\bf W}^2 \|_{2\al', 2m -B}^{1/2} \ge \frac{\rho}{\ve}\})  
\Bigr].
\label{outer.ineq}
\end{eqnarray}
Here $r$ is a sufficiently large integer
and ${\rm Cap}_{2,r}$ is the capacity associated with ${\bf D}_{2,r}$.
Recall that $\theta_{x,a}^{\ve}$
is associated with $ \delta_a (Y^{\ve}_1)$ via Sugita's theorem.

By the large deviation estimate for capacities in \cite{in2},
the second factor on the right hand side of (\ref{outer.ineq})
is known to be  dominated by $\exp (-c (\rho/\ve)^2)$ when $\rho/\ve$
is sufficiently large.
Here, $c = c(\al', m, 2, r)$ is a positive constant.

Suppose  that
\begin{equation}
\| \delta_a (Y^{\ve}_1) \|_{2,-r} = O(\ve^{-\nu})
\qquad
\mbox{ as $\ve \searrow 0$}
\label{pol.dec}
\end{equation}
holds for some $r \in {\mathbb N}$ and some $\nu >0$.
Then, from (\ref{outer.ineq}) and (\ref{pol.dec}), we see that
\begin{eqnarray}
\limsup_{\ve \searrow 0}  \ve^2 \log \mu_{x,a}^{\ve} ( A \cap B^{\prime}_{\rho}({\bf 0})^c)
=
\limsup_{\ve \searrow 0}  \ve^2 \log \mu_{x,a}^{\ve} ( B^{\prime}_{\rho}({\bf 0})^c)
\le
-c \rho^2.
\nn
\end{eqnarray}
Hence, 
\begin{eqnarray}
\limsup_{\ve \searrow 0}  \ve^2 \log \mu_{x,a}^{\ve} ( A)
\le
\Bigl( - \inf
\{
\frac{  \| h\|^2_{{\cal H}}  }{2}
~|~ 
h \in {\cal M}^{x,a}, \, {\cal L} (h) \in A 
\}
\Bigr)
 \vee (-c \rho^2).
\nn
\end{eqnarray}
Letting $\rho \to \infty$, we have 
\begin{eqnarray}
\limsup_{\ve \searrow 0}  \ve^2 \log \mu_{x,a}^{\ve} ( A)
\le
 - \inf
\{
\frac{  \| h\|^2_{{\cal H}}  }{2}
~|~ 
h \in {\cal M}^{x,a}, \, {\cal L} (h) \in A 
\},
\nn
\end{eqnarray}
which is the desired upper estimate.

Now, it remains to prove (\ref{pol.dec}).
We use the integration by parts formula as in (\ref{ups1.eq})--(\ref{ups2.eq}) in Step 1.
However, it is actually easier this time.
(The same symbols are used below.)
Let $Q \in {\bf D}_{\infty}$ be arbitrary.  
We have 
\begin{align}
{\mathbb E}[Q \cdot \delta_a (Y^{\ve}_1) ]
&=
{\mathbb E}[Q  \cdot
( \partial_1^2 \cdots \partial_l^2 )(G \circ \Pi_{{\cal V}} ) (X^{\ve}_1)]
%\nn\\
%&=
=
{\mathbb E}[F_Q^{\ve}  \cdot  (G \circ \Pi_{{\cal V}} ) (X^{\ve}_1) ].
\nn
\end{align}
Here, $F_{Q}^{\ve} (w)$ is a polynomial 
in components of 
(i)~$X^{\ve}_1$ and its derivatives,
(ii)~$Q$ and its derivatives,
(iii)~ $\sigma_{ X^{\ve}_1}$, which is a Malliavin covariance matrix of $X^{\ve}_1$,
and (iv)~ $\gamma_{ X^{\ve}_1}= (\sigma_{ X^{\ve}_1})^{-1}$.
Note that derivatives of $\gamma_{ X^{\ve}_1}$ do not appear.
The right hand side of (\ref{ups2.eq}) 
is not a generalized expectation anymore.

Note that $F_Q$ is linear in $Q$ and there exists $r \in {\mathbb N}$ such that 
the order of derivatives of $Q$ which are involved in the expression of $F_Q^{\ve}$
is bounded from above by $r$.
Combining these with (\ref{ks_malcov.ineq}), we see that 
\[
|{\mathbb E}[Q \cdot \delta_a (Y^{\ve}_1) ] | 
\le 
C \ve^{-\nu} \|Q\|_{2,r}
\qquad
\quad 
(Q \in {\bf D}_{\infty})
\]
for some $\nu >0$ and $C>0$, which are independent of $\ve$ and $Q$.
Since ${\bf D}_{\infty}$ is dense in ${\bf D}_{2,r}$, 
we obtain (\ref{pol.dec}).
This completes the proof of the upper estimate.

%%\newpage
%%%%%%%%%%%%%%%%%%%%%%%%%%%%%%%%%%%%%%%%%%%%%%%%%%%%%%
%%%%%%%%%%%%%%%%%%%%%%%%%%%%%%%
\section{Proof of Corollary \ref{co.pFW} }
%%%%%%%%%%%%%%%%%%%%%%%%%%%%%%%%%%%%%%%%%%%%%%%%%%%%%%%%
%%%%%%%%%%%%%%%%%%%%%%%%%%%%%%%

In this section we prove Corollary \ref{co.pFW}.
Since we are familiar with probability measures, we prove the second assertion,
from which the first assertion immediately follows.
When the vector fields $A_i ~(0 \le i \le d)$ are bounded, 
we can use Lyons continuity theorem for $C_b^3$-vector fields 
and the proof is quite simple due to the contraction principle. 
However, 
when the vector fields  have linear growth, 
we will rely on a cut-off argument and our proof looks a little bit complicated.
(In this section we will write $\lambda^{\ve}_t =\ve^2 t$.)

\subsection{Bounded case}

Consider the product measure $\hat\mu^{\ve}_{x,a} \otimes \delta_{\lambda^{\ve}}$ 
on $G\Omega^H_{\al} ({\mathbb R}^d) \times C_0^{1 -H}([0,1], {\mathbb R})$.
Since the second component $\lambda^{\ve}$ is deterministic 
and continuous in $\ve$,
it follows from (the H\"older version of) Theorem \ref{tm.main}, {\rm (ii)} that 
$\{  \hat\mu^{\ve}_{x,a} \otimes \delta_{\lambda^{\ve}} \}_{0 <\ve \le 1 }$
also satisfies an LDP with a good rate function $J$.
Here, the effective domain of  $J$ is ${\cal L} ( {\cal M}^{x,a} ) \times \{ 0\}$ 
and $J ({\cal L}(h), 0 ) = \|h\|^2_{{\cal H}} /2$
 for $h \in {\cal M}^{x,a}$.

Let $\Phi^{\prime} : G\Omega^H_{\al } ({\mathbb R}^d) \times C_0^{1 -H}([0,1], {\mathbb R})
\to G\Omega^H_{\al} ({\mathbb R}^N)$
be the Lyons-It\^o map 
associated with $A_i ~(0 \le i \le d)$.
Then, 
$
\hat{Z}^{\ve} (\,\cdot\, , z, w) = z+ \Phi^{\prime} (\ve {\cal L} (w) ,  \lambda^{\ve})^1
$
quasi-surely.
(See \cite{in2}.)
In particular, they coincide $\hat\theta^{\ve}_{x,a}$-almost surely.
Then, 
$\tilde{Z}^{\ve}(\,\cdot\,, z)_* [\hat\theta^{\ve}_{x,a}] 
=  (z+ (\Phi^{\prime})^1 )_* [ \hat\mu^{\ve}_{x,a} \otimes \delta_{\lambda^{\ve}}]$.
Note also that $\zeta (h) =  z+ \Phi^{\prime} ( {\cal L} (h) , 0)^1$.
Now, by the contraction principle (Theorem 4.2.1, \cite{dzbk}), 
we can easily show  Corollary \ref{co.pFW}, {\rm (ii) in this case}.

%%%%%%%%%%%%%%%%%%%%%%%%%%%%%%%%%%%%%%%%%%%%%%%%%%%%%%%%
%%%%%%%%%%%%%%%%%%%%%%%%%%%%%%
%  

\subsection{Linearly growing case}

In this case it is not so easy to see whether
the Lyons-It\^o map $\Phi^{\prime}$ is everywhere-defined continuous map or not.
However, as we mentioned before, 
it is well-defined and continuous around (the lift of) Cameron-Martin space. 
Hence, we use a cut-off argument 
and a modified version of the contraction principle (Lemma \ref{lm.cntr} below).
This method has already been used for the usual Freidlin-Wentzell type large deviations
when the coefficient vector fields admit linear growth (see \cite{in5}).

The following lemma is a slight modification of the contraction principle for LDPs
and is formulated in a general setting.
It states that, if the map is continuous around the effective domain of the good rate function, 
then  the contraction principle still holds.
The map need not be continuous everywhere.
\begin{lm}\label{lm.cntr}
Let $S$ and $\hat{S}$ be polish spaces and let $f:S \to \hat{S}$ be a measurable map.
We assume that $\{ \mu_{\ve} \}_{\ve >0}$
is a family of probability measures on $S$ 
which satisfies 
an LDP with a good rate function $J$ as $\ve \searrow 0$.
Let ${\cal D} =\{ a \in S ~|~ J(a) <\infty \}$ be the effective domain of $J$.
Assume further that 
there is an open subset $U$ of $S$ such that ${\cal D} \subset U$
and $f|_U$ is continuous.
Then, $\{ \mu_{\ve} \circ f^{-1} \}_{\ve >0}$ satisfies 
an LDP with a good rate function $\hat{J}$ as $\ve \searrow 0$,
where 
$\hat{J}(b) = \inf\{ J(a) ~|~ a \in f^{-1}(\{b\})  \}$.
\end{lm}

\Proof
For a proof, see Lemma 2.3 \cite{in5}, for instance.
(One can also prove this lemma by hand, since the proof is not so different from the one 
for the standard contraction principle.)
\QED

%%%%%%%%%%%%%%%%%%%%%%%%%%%%
%%%%%%%%%%%%%%%%%%%%%%%%%%%%
%\vspace{5mm}
%%%%%%%%%%%%%%%%%%%%%%%%%%%%%%
%  

Now we discuss a refinement of the Wong-Zakai approximation.
For $k =1,2, \ldots$ and $w \in {\cal W}$, 
$w(k)$ stands for the $k$th dyadic polygonal approximation as before.
We consider the following ODE in the Riemann-Stieltjes sense.
\begin{equation}\label{wz_z.eq}
dz(k)^{\ve}_t = \ve \sum_{i=1}^d  A_i ( z(k)^{\ve}_t)  dw(k)_t^i  + \ve^2  A_0 (z(k)^{\ve}_t)   dt
\qquad
\mbox{with \quad $z(k)^{\ve}_0 =z$}
\end{equation}
It is well-known that if $A_i~(0 \le i \le d)$ satisfies {\bf (A1)},
then for each fixed $\ve$ it holds that
$\lim_{k \to \infty} \sup_{0 \le t \le 1} |Z^{\ve}_t - z(k)^{\ve}_t| = 0$
for $\mu$-a.a. $w$.
However, this convergence actually takes place quasi-surely.

In the next lemma 
$O$ denotes the largest open subset of $G\Omega^H_{\al} ({\mathbb R}^d) \times C_0^{1 -H}([0,1], {\mathbb R})$ 
on which $\Phi^{\prime}$ is well-defined.
As we already explained, 
$\Phi^{\prime}$ is continuous from $O$ to $G\Omega^H_{\al} ({\mathbb R}^N)$.

%%%%%%%%%%
%\vspace{5mm}
%%%%%%%%%%

\begin{lm}\label{lm.qs_cinc} 
Assume that 
the vector field $A_i$ satisfies {\bf (A1)} for any $0 \le i \le d$
and fix any $\ve \in (0,1]$.
Then, quasi-surely, 
\[
\lim_{k \to \infty} \sup_{0 \le t \le 1} |\tilde{Z}^{\ve}_t - z(k)^{\ve}_t| = 0.
\]
Here, $\tilde{Z}^{\ve}=\tilde{Z}^{\ve} (\,\cdot\, , z, \,\cdot\,)$ 
stands for the $\infty$-quasi-sure modification of 
$Z^{\ve}$ defined in (\ref{tildeZ.def}).
Moreover, we have 
\[
\tilde{Z}^{\ve} (\,\cdot\, , z, w)= z + \Phi^{\prime} (\ve {\bf W}, \lambda^{\ve})^1
\qquad
\mbox{ quasi-surely on $\{ w \in {\cal W}~|~ (\ve {\bf W}, \lambda^{\ve}) \in O\}$.}
\]
\end{lm}

\Proof
For any $0<a<b$, choose 
a smooth, non-increasing function $\chi_{a,b}:{\mathbb R} \to [0,1]$
such  that 
$\chi_{a,b} =1$ on $(-\infty, a]$,  $\chi_{a,b} >0$ on $(-\infty, b)$,
and $\chi_{a,b} =1$ on $[b, \infty)$.
For $\nu=1,2,\ldots$, we set 
$A_i^{\nu} (z) = \chi_{\nu +1, \nu +2}( |z| ) A_i (z)$ for $z \in {\mathbb R}^N$ and $0 \le i \le d$.
Clearly, 
$A_i^{\nu}$ is of $C^{\infty}_b$ and agrees with $A_i$ on $\{ z~|~ |z| \le \nu +1\}$.
Consider the scaled SDE (\ref{sc_sde_A.def}) and 
its approximating ODE (\ref{wz_z.eq})
with their  coefficient vector fields being replaced by $A_i^{\nu}$.
The solutions are denoted by $Z^{\ve, \nu}$ and $z^{\ve, \nu}(k)$, respectively.
The Lyons-It\^o map 
associated with the new coefficients is denoted by $\Phi^{\prime}_{\nu}$,
which is defined everywhere and 
continuous from of $G\Omega^H_{\al} ({\mathbb R}^d) \times C_0^{1 -H}([0,1], {\mathbb R})$ 
to 
$G\Omega^H_{\al} ({\mathbb R}^N)$.
Then, 
$\tilde{Z}^{\ve, \nu} = z + \Phi^{\prime}_{\nu}(\ve {\bf W}, \lambda^{\ve})^1$, quasi-surely.

Take any $\nu$ such that $\nu \ge |z|$ and we will denote the sup-norm by $\|  \,\cdot\,\|_{\infty} $.
By a standard argument for stopping times, 
$\{ w~|~  \|  Z^{\ve}\|_{\infty} <\nu \}=\{ w~|~  \|  Z^{\ve, \nu}\|_{\infty} <\nu \}$, $\mu$-a.s.
and 
$Z^{\ve} = Z^{\ve, \nu}$, $\mu$-a.s. on this subset.
Hence, 
for any $0<a<b <\nu +1$,
$\chi_{a,b}( \|  Z^{\ve}\|_{\infty} ) Z^{\ve} 
= \chi_{a,b}( \|  Z^{\ve, \nu}\|_{\infty} ) Z^{\ve, \nu}$, $\mu$-a.s. on ${\cal W}$.
By the uniqueness of quasi-continuous modification, 
$\chi_{a,b}( \|  \tilde{Z}^{\ve}\|_{\infty} ) \tilde{Z}^{\ve} 
= \chi_{a,b}( \|  \tilde{Z}^{\ve, \nu}\|_{\infty})   \tilde{Z}^{\ve, \nu}$, quasi-surely on ${\cal W}$.
Assume that $z\neq 0$, since the case $z=0$ can be shown with trivial modification.
Since $ \tilde{Z}^{\ve}$ and 
$ \tilde{Z}^{\ve, \nu}$
can never be a zero path, this implies that 
$\{ w~|~  \|  Z^{\ve}\|_{\infty} <b \}=\{ w~|~  \|  Z^{\ve, \nu}\|_{\infty} <b \}$, quasi-surely,
for any $b \in (0, \nu +1)$.
Taking $b=a$ and 
using the above equality once again, we have 
$ Z^{\ve}=  Z^{\ve, \nu}$
quasi-surely on the above subset for any $b \in (0, \nu +1)$.
(Below, we will use this fact with $b =\nu$.)

If $w \in \{ w~|~  \|  Z^{\ve}\|_{\infty} <\nu \}=\{ w~|~  \|  Z^{\ve, \nu}\|_{\infty} <\nu \}$
and 
admits a rough path lift with respect to $\alpha$-H\"older rough path topology,
then it is easy to see that
\begin{align}
Z^{\ve} (\,\cdot\, , z,w ) = Z^{\ve, \nu} (\,\cdot\, , z,w ) 
&= z + \Phi^{\prime}_{\nu}(\ve {\cal L}(w), \lambda^{\ve})^1
%\nn\\&
=
\lim_{m \to \infty}  \Bigl( 
z + \Phi^{\prime}_{\nu}(\ve {\cal L}( w(k)), \lambda^{\ve})^1
\Bigr),
%=
%\lim_{m \to \infty} z^{\ve, \nu}(m)(\,\cdot\, , z,w ).
\nn
\end{align}
where we have used Lyons' continuity theorem for $\Phi^{\prime}_{\nu}$.
Since 
$z + \Phi^{\prime}_{\nu}(\ve {\cal L}( w(k)), \lambda^{\ve})^1 = z^{\ve, \nu}(k)(\,\cdot\, , z,w )$ 
stays inside the ball of radius $\nu +1$ for sufficiently large $k$,
it holds that $z^{\ve, \nu}(k)(\,\cdot\, , z,w ) = z^{\ve}(k)(\,\cdot\, , z,w )$.
Thus, we have shown the refinement of the Wong-Zakai approximation
on the set $\{ w~|~  \|  Z^{\ve}\|_{\infty} <\nu \}$ and, by 
taking the union with respect to $\nu$, on the whole Wiener space, too. 

The proof of 
the second assertion of the lemma is quite similar.
We just need to note that 
$\Phi^{\prime}$ is continuous on the open set $O$  and 
that $z + \Phi^{\prime}_{\nu}(\ve {\cal L}(w), \lambda^{\ve})^1 
= z + \Phi^{\prime} (\ve {\cal L}(w), \lambda^{\ve})^1$
as long as 
it stays inside the ball of radius $\nu +1$ for sufficiently large $k$.
\QED

\noindent
{\it Proof of the linear growth case of Corollary \ref{co.pFW}.}~
For simplicity of notation we prove the case $z=0$ only.
We extend $(\Phi^{\prime})^1: O \to C^{\alpha -H}([0,1], {\mathbb R}^N)$  
by setting 
$\Phi^{\prime}({\bf w}, \lambda)^1 = 0$ if $({\bf w}, \lambda) \notin O$.
Note that 
$O$ contains ${\cal L}({\cal H}) \times \{0\}$, which in turn contains 
the effective domain of the rate function $\hat{I}_1$ in Theorem \ref{tm.main}.
Then, 
by Lemma \ref{lm.cntr}, 
the push-forward measure of $\hat{\mu}_{x,a}^{\ve}$
by the map $(\Phi^{\prime})^1$
satisfies an LDP with a good rate function $\hat{I}_2$.

Fix $\ve \in (0,1]$.
On the probability space $({\cal W},\hat{\theta}_{x,a}^{\ve})$, 
we have two $C^{\alpha -H}([0,1], {\mathbb R}^N)$-valued random maps.
One is $\tilde{Z}^{\ve}$.
The other is 
$w \mapsto \Phi^{\prime}(\ve {\bf W}, \lambda^{\ve})^1$.
The push-forward measure of $\hat{\theta}_{x,a}^{\ve}$
by the latter map is $(\Phi^{\prime})^1_* [\hat{\mu}_{x,a}^{\ve}]$, 
which we have just discussed.

Let us consider the set on which these two maps disagree.
By Lemma \ref{lm.qs_cinc} and the fact that $\hat{\theta}_{x,a}^{\ve}$ does not 
charge a slim set,
\begin{eqnarray*}
\hat{\theta}_{x,a}^{\ve} (\{ w~|~ \tilde{Z}^{\ve} (\,\cdot\, ,z,w)
\neq \Phi^{\prime}(\ve {\bf W}, \lambda^{\ve})^1\})
&\le& 
\hat{\theta}_{x,a}^{\ve} (\{ w~|~ (\ve {\bf W}, \lambda^{\ve}) \notin O\})
\nn\\
&=& 
  \hat\mu^{\ve}_{x,a} \otimes \delta_{\lambda^{\ve}} (O^c).
\end{eqnarray*}
Our main theorem (Theorem \ref{tm.main}) implies that $\limsup_{\ve \searrow 0} 
\ve^2 \log \hat\mu^{\ve}_{x,a} \otimes \delta_{\lambda^{\ve}} (O^c) = -\infty$.
Therefore, these two random maps are exponentially equivalent in the sense of 
Definition 4.2.10, \cite{dzbk}. 
We see from 
Theorem 4.2.13, \cite{dzbk} 
that 
$(\tilde{Z}^{\ve})_* [\hat{\theta}_{x,a}^{\ve}]$ also satisfies 
an LDP with the same good rate function $\hat{I}_2$.
This completes the proof of Corollary \ref{co.pFW}, {\rm (ii)}.
\toy

%\newpage
%%%%%%%%%%%%%%%%%%%%%%%%%%%%%%%%%%%%%%%%%%%%%%%%%%%%%%%%
%%%%%%%%%%%%%%%%%%%%%%%%%%%%%%
%  references
%%%%%%%%%%%%%%%%%%%%%%%%%%%%%%%%%%%%%%%%%%%%%%%%%%%%%%%%
%%%%%%%%%%%%%%%%%%%%%%%%%%%%%%%

\vspace{15mm}

\begin{flushleft}
\begin{tabular}{ll}
Yuzuru INAHAMA
\\
Graduate School of Mathematics,   Nagoya University,
\\
Furocho, Chikusa-ku, Nagoya 464-8602, JAPAN.
\\
Email: {\tt inahama@math.nagoya-u.ac.jp}
\end{tabular}
\end{flushleft}

\end{document}